\documentclass{article}
\usepackage{amssymb}
\usepackage{amsmath}
\usepackage[dvips]{graphicx}
\numberwithin{equation}{section}
\newcommand{\RR}{{\mathbb{R}}}
\newcommand{\NN}{{\mathbb{N}}}

\newcommand{\meanint}{{\int{\mkern-19mu}-}}
\def\ring{\mathaccent"0017 }

\title{\bf\Large  Recent progress in  elliptic equations and systems of arbitrary  order with rough coefficients in Lipschitz  domains}

{ \date\ }

\newtheorem{lemma}{Lemma}
\newtheorem{theorem}{Theorem}
\newtheorem{proposition}{Proposition}
\newtheorem{corollary}{Corollary}

\begin{document}
\maketitle
\vspace{-2cm}
\centerline{\scshape Vladimir Maz'ya}
\medskip
{\footnotesize
 \centerline{Department of Mathematical Sciences}
   \centerline{University of Liverpool, M$\&$O Building}
   \centerline{Liverpool, L69 3BX, UK}
\medskip
 \centerline{Department of Mathematics}
   \centerline{Link\"oping University}
   \centerline{SE-58183 Link\"oping, Sweden}}
   \medskip
   \centerline{\scshape Tatyana Shaposhnikova}
\medskip
{\footnotesize
 \centerline{Department of  Mathematics}
   \centerline{Link\"oping University}
   \centerline{SE-58183 Link\"oping, Sweden}} 
\bigskip
\bigskip

{\bf Abstract.}
This is a survey of  results mostly relating   elliptic equations and systems of arbitrary even order with rough coefficients in Lipschitz graph  domains. Asymptotic properties of solutions at a point of a Lipschitz boundary are also discussed.

\section*{Introduction}

The fundamental role 
in the theory of linear elliptic equations and systems is played by results on regularity of solutions up to  the boundary of a domain. The following classical example serves as an illustration.

Consider the Dirichlet problem

$$ \left\{\begin{array}{l}
\Delta u=0\hskip 0.20in\mbox{in}\quad \Omega,\\[.25cm]
{\rm Tr} \, u=f \hskip 0.28in\mbox{on}\quad
\partial\Omega, 
\hskip 0.20in
\end{array}\right.
$$
where $\Omega$ is a bounded domain with smooth boundary and ${\rm Tr}\, u$ stands for the boundary value (trace) of $u$. 
Let $u\in L_p^1(\Omega)$, $1<p<\infty$, that is
$$\int_{\Omega} |\nabla u|^p dx <\infty$$
and let $f$ belong to the Besov space  $B_p^{1-1/p}(\partial\Omega)$ with the seminorm
$$\Bigl(\int_{\partial\Omega}\int_{\partial\Omega} \frac{|f(x) - f(y)|^p}{|x-y|^{n+p-2}} d\sigma_xd\sigma_y \Bigr)^{1/p}. $$
It is well-known that
$${\rm Tr}\, L_p^1(\Omega) = B_p^{1-1/p}(\partial\Omega).$$
Moreover, the harmonic extension  $u$ of $f\in B_p^{1-1/p}(\partial\Omega)$ belongs to $L_p^1(\Omega)$ and the norm of the gradient $\nabla u$ in $L_p(\Omega)$ is equivalent to the above seminorm in $B_p^{1-1/p}(\partial\Omega)$.

\smallskip

This fact highlights the following topics
 of interest:

\begin{itemize}
\item Replace the domain with smooth boundary by a domain in more general  class and study the effect of irregularities of $\partial\Omega$.

\item Replace $\Delta$ by a more general  elliptic 
 operator with variable coefficients, and study the impact of low regularity 
assumptions on the coefficients.

\item Understand the correlation between the smoothness of data and 
the smoothness of solutions.

\end{itemize}

Another theme, somehow related to these, is: 

\begin{itemize}
\item Describe  the local behaviour of solutions near a  boundary or  interior
point of the domain $\Omega$.
\end{itemize}

\smallskip

In the present paper we survey   results in the directions just mentioned. Most of them were obtained during the last decade and concern  elliptic equations and  systems of arbitrary  order. Special attention is paid to the Stokes system. 
The selection of topics is partly influenced by our involvement in their study. Here is the plan of the article.

\smallskip

Section 1 is dedicated to  weak solutions with Besov boundary data, with coefficients of differential operators and the unit normal to the boundary in  classes close to $BMO$. Results of a similar nature for the Stokes system are discussed in Section 2. In Section 3, strong solutions  in Sobolev spaces are considered. Here sharp  additional conditions on the Lipschitz boundary are reviewed. In particular, in Subsection 3.2 we speak about strong solvability of the Stokes system. 
In Section 4, asymptotic formulas for solutions of the Dirichlet problem near an isolated point of the Lipschitz boundary and at a point in the domain are dealt with. 

\section{Weak solvability of higher order elliptic systems with coefficients close to $BMO$ in Lipschitz domains}

\subsection{Background}
 The present section is mostly based  on results of the paper by V. Maz'ya, M. Mitrea and T. Shaposhnikova \cite{[MMS1]}. We start with mentioning earlier works.

\smallskip

 The  basic case of the  Laplacian in arbitrary Lipschitz graph domains in 
$\RR^n$, is treated in the work of B.\,Dahlberg and 
C.\,Kenig {\cite{[DK]}}, in the case of $L_p$-data, and D.\,Jerison and 
C.\,Kenig {\cite{[JK]}}, in the case of data in the Besov space $B_p^s$ with $0<s<1$. The local regularity in the Sobolev class $W_p^2$ of solutions to second order equations with coefficients in $VMO\cap L_\infty$ was established by F. Chiarenza, M. Frasca, P. Longo \cite{[CFL]}. 

 \smallskip
 
 In spite of substantial progress in recent years, there remain many basic 
open questions  for higher order equations, even in the case of 
{constant coefficient} operators in Lipschitz domains.  
One significant problem is to determine the sharp range of $p$'s for which 
the Dirichlet problem for strongly elliptic systems with $L_p$-boundary data 
is well-posed. In {\cite{[PV]}}, 
J.\,Pipher and G.\,Verchota have developed a $L_p$-theory for real, constant 
coefficient, higher order systems 
$$L=\sum_{|\alpha|=2m}A_\alpha D^\alpha$$ 
when 
$p$ is near $2$, i.e., $2-\varepsilon<p<2+\varepsilon$ with $\varepsilon>0$ 
depending on the Lipschitz character of $\Omega$. On p.\,2 of {\cite{[PV]}} 
the authors ask whether the $L_p$-Dirichlet problem for these operators is 
solvable in a given Lipschitz domain for 
$p\in (2-\varepsilon,\frac{2(n-1)}{n-3}+\varepsilon)$, and a positive answer 
has been  given by Z.\,Shen in {\cite{[Sh]}}. Let us also mention 
 the work {\cite{[AP]}} of V.\,Adolfsson and J.\,Pipher who have dealt 
with the Dirichlet problem for the biharmonic operator in arbitrary Lipschitz 
domains and with data in Besov spaces, as well as {\cite{[Ve]}} where G.\,Verchota 
formulates and solves a Neumann-type problem for the bi-Laplacian in Lipschitz
domains and with boundary data in $L_2$. In  {\cite{[MMT]}} D.\,Mitrea, M.\,Mitrea,  and M.\,Taylor, 
treat the Dirichlet problem for strongly elliptic systems 
of second order in an arbitrary Lipschitz subdomain $\Omega$ of a (smooth) 
Riemannian manifold and with boundary data in $B^s_p(\partial\Omega)$, 
when $2-\varepsilon<p<2+\varepsilon$ and $0<s<1$.

\smallskip

We mention some recent results on the second order elliptic equations and systems with coefficients in $VMO$ 
due to 
 G. Di Fazio \cite{[DF]}, L. Caffarelli and I. Peral \cite{[CP]}, B. Stroffolini \cite{[St]},  D. Guidetti  \cite{[Gu]}, for higher order equations see also D. Palagachev and L. Softova \cite{[PS]}.

\subsection{Domains and function spaces}
 Let us turn to the article \cite{[MMS1]}. 
 We make no notational distinction between 
 spaces of scalar-valued functions and their natural 
counterparts for vector-valued functions.

\smallskip

Recall that a domain $\Omega$ is called  Lipschitz  graph if its 
boundary can be locally described by means of (appropriately rotated
and translated) graphs of real-valued Lipschitz functions.  

\smallskip

It is shown by S. Hofmann, M.  Mitrea,  and M. Taylor in \cite{[HMT]} that $\Omega$ is a  Lipschitz  graph domain if and only if it has finite perimeter in the sense of De Giorgi (see \cite{[DG1]}, \cite {[Fe]}, \cite{[BM]}) and (i) there are continuous (or, equivalently, smooth) vector fields that are transversal to the boundary and (ii) the necessary  condition $\partial\Omega = \partial\overline\Omega$ is fulfilled. 

\smallskip

By \cite{[HMT]}, a bounded nonempty domain of finite perimeter for which $\partial\Omega = \partial\overline\Omega$ is Lipschitz graph if and only if 
$$\inf \{\|\nu-\omega\|_{L_\infty(\partial\Omega)}: \,\, \omega =(\omega_1, \ldots , \omega_n)\in C^0(\partial\Omega), \,\,\,  |\omega| =1\, \,\,  {\rm on} \, \, \, \partial\Omega\} <\sqrt{2}$$
with $\nu$ being the outward normal to $\partial\Omega$.

\smallskip

Everywhere in this section we assume that $\Omega$ is a Lipschitz graph  domain  in ${\Bbb{R}}^n$, 
with compact closure $\overline{\Omega}$ and with 
outward unit normal  $\mathbf{\nu}=(\nu_1,...,\nu_n)$.  Let $m$ be an integer.  Consider 
the   operator 
\begin{equation}\label{star}
{\mathcal L}(X,D_X)\,{\mathcal U}
:=\sum_{|\alpha|=|\beta|=m}D^\alpha(A_{\alpha\beta}(X)D^\beta{\mathcal U}), \quad X\in \Omega,
\end{equation}
with the data
\begin{equation}\label{eq1}
\displaystyle{\frac{\partial^k{\mathcal U}}{\partial\nu^k}}=g_k
\,\,\quad\mbox{on}\,\,\partial\Omega,\qquad 0\leq k\leq m-1.
\end{equation}

\noindent  The coefficients $A_{\alpha\beta}$ are square  matrices
with measurable, complex-valued entries, for which $\exists\,\kappa>0$ 
\begin{equation}\label{A-bdd}
\sum_{|\alpha|=|\beta|=m}\|A_{\alpha\beta}\|_{L_\infty(\Omega)}\leq\kappa^{-1}
\end{equation}
\noindent  and such that the coercivity condition
$$ {\rm Re}\,\int_\Omega\sum_{|\alpha|=|\beta|=m}\langle A_{\alpha\beta}(X) 
D^\beta\,{ V}(X),\,D^\alpha\,{ V}(X)\rangle\,dX 
 \geq\kappa\sum_{|\alpha|=m}\|D^\alpha\,{V}\|^2_{L_2(\Omega)}$$
\noindent  holds for all  complex vector-valued 
functions ${ V}\in C^\infty_0(\Omega)$.

\smallskip

Let $\mathcal U$ belong to the usual Sobolev space $W_p^m(\Omega)$. 
It is natural to take 
\begin{eqnarray*}
\frac{\partial^k{\mathcal U}}{\partial\nu^k}
:=\sum_{|\alpha|=k}\frac{k!}{\alpha!}\,
\nu^\alpha\,{\rm Tr}\,[D^\alpha{\mathcal U}],\quad 0\leq k\leq m-1,
\end{eqnarray*}
\noindent  where 
$\nu^\alpha:=\nu_1^{\alpha_1}\cdots\nu_n^{\alpha_n}$  if 
$\alpha=(\alpha_1,...,\alpha_n)$.
Now, if  $p\in(1,\infty)$,  $a\in(-1/p,1-1/p)$
are fixed and 
$$
\rho(X):={\rm dist}\,(X,\partial\Omega)
$$
the space  $W^{m,a}_p(\Omega)$ is 
defined as the space of vector-valued functions for which 
\begin{equation}\label{1.6}
\sum_{|\alpha|\leq m}\int_\Omega|D^\alpha{\mathcal U}(X)|^p 
\rho(X)^{pa}\,dX<\infty.
\end{equation}
Further we set
\begin{equation}\label{w}
V^{m,a}_p(\Omega):=\mbox{the closure of }C^\infty_0(\Omega)\mbox{ in }
W^{m,a}_p(\Omega)
\end{equation}
\noindent  
and introduce the dual space
\begin{eqnarray*}
V^{-m,a}_p(\Omega):=\bigl(V^{m,-a}_{p'}(\Omega)\bigr)^*
\end{eqnarray*}

\smallskip

For any ${\mathcal U}\in W^{m,a}_p(\Omega)$, 
the  traces of $D^\alpha{\mathcal U}$, $0\leq |\alpha|\leq m-1$, 
exist in 
$B_p^s(\partial\Omega)$, where  $s:=1-a-1/p$, $0< s<1$ (see E. Gagliardo \cite{[Gag]} for $a=0$ and S. Uspenski\u{\i} \cite{[Usp]}).
Recall that $f\in L_p(\partial\Omega)$ belongs to $B_p^s(\partial\Omega)$ if and only if
\begin{equation}\label{eq2}
\int_{\partial\Omega}\int_{\partial\Omega}
\frac{|f(X)-f(Y)|^p}{|X-Y|^{n-1+sp}}\,d\sigma_Xd\sigma_Y<\infty.
\end{equation}

\smallskip

The above definition takes advantage of the Lipschitz manifold structure 
of $\partial\Omega$ which allows one to define smoothness spaces of index  less than 
 $1$. 
{The approach is no longer effective when the order of smoothness exceeds
 $1$}. 
 
 \smallskip
 
 Let us illustrate the necessity of working with boundary data 
 different from those in spaces of traces of usual Sobolev spaces   by considering 
Dirichlet problem for the biharmonic operator 
\begin{eqnarray*}
&& {\cal U}\in W_2^2(\Omega),\qquad 
\Delta ^2\,{\cal U}=0\,\,\mbox{in}\,\,\Omega,
\\[4pt]
&& {\rm Tr} \,{\cal U}=g_0\,\,\mbox{on}\,\,\partial\Omega,\quad
\langle\nu,{\rm Tr}\,[\nabla{\cal U}]\rangle 
=g_1\,\,\mbox{on}\,\,\partial\Omega.
\end{eqnarray*}
\noindent  One might be tempted to believe that a natural 
class of boundary data
is  $B_{2}^{3/2}(\partial\Omega)\times B_{2}^{1/2}(\partial\Omega)$, 
 where  by definition  
$B_{2}^{3/2}(\partial\Omega)$ and $B_{2}^{1/2}(\partial\Omega)$ 
are the spaces of traces of functions in $ W_2^2(\Omega)$ and 
$W_2^1(\Omega)$, respectively. 

\smallskip

However, this formulation has  serious drawbacks. 
The first one is that the mapping 
\begin{eqnarray*}
W_2^2(\Omega)\ni {\cal U}\mapsto
\langle\nu,{\rm Tr}\,[\nabla{\mathcal U}]\rangle\in 
B_{2}^{1/2}(\partial\Omega)
\end{eqnarray*}
\noindent  is generally unbounded. In fact, its continuity implies 
 $\nu\in B_{2}^{1/2}(\partial\Omega)$  
which is not necessarily the case for a Lipschitz domain, even for 
the square $S=[0,1]^2$.

\smallskip

Second, this problem may fail to have a solution when 
 $(g_0, g_1)$ is an arbitrary pair in $B_{2}^{3/2}(\partial\Omega)\times 
B_{2}^{1/2}(\partial\Omega)$ 
. Indeed, consider the case  $\Omega=S$  and 
 $g_0=0$, $g_1=1$. It is standard that the main term of the 
asymptotics near the origin of any solution ${\mathcal U}$ in $W_2^2(S)$ 
is given in polar coordinates $(r,\theta)$ by 
\begin{eqnarray*}
\frac{2r}{\pi+2}\left((\theta-\frac{\pi}{2})\sin\theta
-\theta\cos\theta\right).
\end{eqnarray*}
\noindent  Since this function does not belong to $W_2^2(S)$, there 
is no solution in this space.

\smallskip

A new point of view has been 
introduced by H.\,Whitney in {\cite{[Wh]}} who  considered higher order 
Lipschitz spaces on arbitrary closed sets. An extension of this 
circle of ideas pertaining to the full scale of Besov and Sobolev spaces on 
irregular subsets of $\RR^n$ can be found in the book {\cite{[JW]}} by 
A.\,Jonsson and H.\,Wallin. 
 The authors of \cite{[MMS1]} further refined this theory in the context 
of Lipschitz domains. The description of  higher order Besov spaces on the boundary of 
a Lipschitz domain $\Omega\subset\RR^n$ in \cite{[MMS1]} runs as follows.

\smallskip

For $m\in{\mathbb{N}}$, $p\in(1,\infty)$, $s\in(0,1)$, the  
space $\dot{B}^{m-1+s}_{p}(\partial\Omega)$ 
is introduced as the closure of 
\begin{eqnarray*}
\Bigl\{(D^\alpha\,{\mathcal V}|_{\partial\Omega})
_{|\alpha|\leq m-1}:\,{\mathcal V}\in C^\infty_0(\RR^n)\Bigr\}
\end{eqnarray*}
\noindent 
 in $B_{p}^s(\partial\Omega)$. An equivalent characterization 
of $\dot{B}^{m-1+s}_p(\partial\Omega)$ which involves higher order Taylor remainder  in place of
$f(X)-f(Y)$ in (\ref{eq2})
   runs as follows  (see  Sect. {7.1} of \cite{[MMS1]}).

\smallskip

For a collection of 
families $\dot{f}=\{f_\alpha\}_{|\alpha|\leq m-1}$ of 
measurable functions defined on $\partial\Omega$,   there is the equivalence relation
\begin{eqnarray}\label{Bes-Nr}
\|\dot{f}\|_{\dot{B}^{m-1+s}_p(\partial\Omega)}
& \sim & \sum_{|\alpha|\leq m-1}\|f_\alpha\|_{L_p(\partial\Omega)}
\\[6pt]
&& +\sum_{|\alpha|\leq m-1}\Bigl(\int_{\partial\Omega}\int_{\partial\Omega}
\frac{|R_\alpha(X,Y)|^p}{|X-Y|^{p(m-1+s-|\alpha|)+n-1}}\,d\sigma_Xd\sigma_Y
\Bigr)^{1/p}<\infty,
\nonumber
\end{eqnarray}
where
\begin{equation}\label{reminder}
R_\alpha(X,Y):=f_\alpha(X)-\sum_{|\beta|\leq m-1-|\alpha|}\frac{1}{\beta!}\,
f_{\alpha+\beta}(Y)\,(X-Y)^\beta,\qquad X,Y\in\partial\Omega,
\end{equation}

\noindent It is standard to prove that $\dot{B}^{m-1+s}_p(\partial\Omega)$ is a 
Banach space. Also, trivially, for any constant $\kappa>0$, 
\begin{equation}\label{Bes-Nr-eq}
\sum_{|\alpha|\leq m-1}\|f_\alpha\|_{L_p(\partial\Omega)}
+\sum_{|\alpha|\leq m-1}\,\,\,\Bigl(\,\,\,\,\,\,\,\,\,
\int\!\!\!\!\!\!\!\!\!\!\!\!\!\!\!\!\!\!\!
\int\limits_{{X,Y\in\partial\Omega}\atop{|X-Y|<\kappa}}
\frac{|R_\alpha(X,Y)|^p}{|X-Y|^{p(m-1+s-|\alpha|)+n-1}}\,d\sigma_Xd\sigma_Y
\Bigr)^{1/p}
\end{equation}
\noindent is an equivalent norm on $\dot{B}^{m-1+s}_p(\partial\Omega)$.

\smallskip

In order to formulate the trace and extension theorem for the spaces $\dot{B}^{m-1+s}_p(\partial\Omega)$, we first give its analogue for lower smoothness spaces which is essentially due to S.\,Uspenski\u{\i} 
 {\cite{[Usp]}}.

\smallskip

\begin{lemma}\label{trace-1}
For each $1<p<\infty$, $-1/p<a<1-1/p$ and $s:=1-a-1/p$, the trace operator
\begin{equation}\label{TR-1}
{\rm Tr}:W^{1,a}_p(\Omega)\longrightarrow B^s_p(\partial\Omega)
\end{equation}
\noindent is well-defined, linear, bounded, onto and has $V^{1,a}_p(\Omega)$ 
as its null-space. Furthermore, there exists a linear, continuous mapping
\begin{equation}\label{Extension}
{\mathcal E}:B^s_p(\partial\Omega)\longrightarrow W^{1,a}_p(\Omega),
\end{equation}
\noindent called extension operator, such that ${\rm Tr}\circ{\mathcal E}=I$
(i.e., a bounded, linear right-inverse of trace). 
\end{lemma}

\smallskip

For higher smoothness see the following assertion which is 
Proposition 7.3 in \cite{[MMS1]}.

\begin{proposition}\label{trace-2}
For $1<p<\infty$, $-1/p<a<1-1/p$, $s:=1-a-1/p\in(0,1)$ and $m\in\NN$, 
define the  higher order  trace operator
\begin{equation}\label{TR-11}
{\rm tr}_{m-1}:W^{m,a}_p(\Omega)\longrightarrow
\dot{B}^{m-1+s}_p(\partial\Omega)
\end{equation}
\noindent by setting 
\begin{equation}\label{Tr-DDD}
{\rm tr}_{m-1}\,\,{\mathcal U}
:=\Bigl\{i^{|\alpha|}\,{\rm Tr}\,[D^\alpha\,{\mathcal U}]\Bigr\}
_{|\alpha|\leq m-1},
\end{equation}
\noindent where the traces in the right-hand side are taken in the sense
of Lemma~\ref{trace-1}. Then {\rm (\ref{TR-11})} is 
a well-defined, linear, bounded operator, which is onto and has 
$V^{m,a}_p(\Omega)$ as its null-space. Moreover, it has a bounded, linear 
right-inverse, i.e., there exists a linear, continuous operator 
\begin{equation}\label{Ext-222}
{\mathcal E}:\dot{B}^{m-1+s}_p(\partial\Omega)
\longrightarrow W^{m,a}_p(\Omega)
\end{equation}
\noindent such that 
\begin{equation}\label{Ext-333}
\dot{f}=\{f_\alpha\}_{|\alpha|\leq m-1}\in\dot{B}^{m-1+s}_p(\partial\Omega)
\Longrightarrow i^{|\alpha|}\,{\rm Tr}\,[D^\alpha({\mathcal E}\,\dot{f})]
=f_\alpha . 
\end{equation}
\end{proposition}

\smallskip

Now,  a necessary condition for the boundary data 
$\{g_k\}_{0\leq k\leq m-1}$  in the Dirichlet problem (\ref{eq1}) 
is that 
\begin{eqnarray}\label{1.9}
\begin{array}{l}
\displaystyle{
\exists\,\dot{f}=\{f_\alpha\}_{|\alpha|\leq m-1}
\in\dot{B}^{m-1+s}_p(\partial\Omega) \quad {\rm such}\,\,{\rm that}}\\[10pt]
\displaystyle{
g_k=\sum_{|\alpha|=k}\frac{k!}{\alpha!}\,\nu^\alpha\,f_\alpha
\quad\mbox{for each}\,\,\,\,0\leq k\leq m-1.}
\end{array}
\end{eqnarray}
\noindent This family $\{g_k\}$ is organized as a Banach space, and is 
denoted by $\dot{W}^{m-1+s}_p(\partial\Omega)$. 

\smallskip

The space  takes a particularly simple form when $m=2$. To describe it we need 
the notation for the tangential derivative
  $\partial/\partial\tau_{jk}$
given by 
\begin{equation}\label{tang-tau}
\frac{\partial}{\partial\tau_{jk}}:=\nu_j\frac{\partial}{\partial x_k}
-\nu_k\frac{\partial}{\partial x_j},\qquad 1\leq j,k\leq n
\end{equation}
and the tangential gradient on the surface $\partial\Omega$
$$\nabla_{\rm tan}:=(\sum_j\nu_j\partial/\partial\tau_{jk})_{1\leq k\leq n}.$$
Then, 
for each Lipschitz graph domain $\Omega\subset\RR^n$ and each $1<p<\infty$, 
$s\in(0,1)$, 
\begin{equation}\label{W2-az}
\dot{W}^{1+s}_p(\partial\Omega)
=\{(g_0,g_1)\in L^1_p(\partial\Omega)\oplus L_p(\partial\Omega):\,
\nu g_1+\nabla_{\rm tan}\,g_0\in B^s_p(\partial\Omega)\}.
\end{equation}

\noindent This has been conjectured to hold (when $s=1-1/p$) by 
A.\,Buffa and G.\,Geymonat on p.\,703 of {{\cite{[BG]}}}.

\subsection{Formulation of the Dirichlet problem}

Broadly speaking, there are two types of 
questions pertaining to the well-posedness of 
the Dirichlet problem in a Lipschitz domain $\Omega$ 
for a divergence form strongly  elliptic system of order 
$2m$ with boundary data in $\dot{W}^{m-1+s}_p(\partial\Omega)$.

\vskip 0.08in
\noindent { { Question I.}} 
Granted that the coefficients of ${\mathcal L}$ 
exhibit a certain amount of smoothness, identify the indices $p$, $s$  
for which this boundary value problem is well-posed. 

\vskip 0.08in
\noindent  {{Question II.}} 
Alternatively, having fixed the indices $s$ 
and $p$, characterize the smoothness of $\partial\Omega$ and 
of the coefficients of ${\mathcal L}$ for which the aforementioned problem 
is well-posed. 

\smallskip

Both questions are discussed in \cite{[MMS1]}.

\smallskip

 B. Dahlberg (see \cite{[Dah1]} has constructed a domain $\Omega\subset\RR^n$, 
bounded with $C^1$-boundary, and a function $f\in C^\infty(\bar{\Omega})$ 
such that, for each $p\in (1, \infty)$, 
$$
\Delta u=f,\,\,u\in W^1_{2}(\Omega),\,\,u|_{\partial\Omega}=0 \quad \Longrightarrow
\partial_j \partial_k u \notin L_p(\Omega),
$$
where $\partial_j = \partial / \partial x_j$. 
In other words, Dahlberg's domain $\Omega$ has the property that
for each $p\in (1, \infty)$, the {\it Poisson problem for the Laplacian 
with homogeneous Dirichlet boundary conditions} fails to be well posed, 
in the sense that it is no longer reasonable 
to expect solutions with  two  derivatives in $L_p(\Omega)$. 

A fundamental issue raised by this counterexample:

{{\it Identify those Sobolev-Besov spaces within which the natural 
correlation between the smoothness of the data and that of the solutions of the Dirichlet problem 
is preserved when the domain in question is allowed to have a minimally 
smooth boundary}}.

\medskip

Consider the Dirichlet problem  for the operator (\ref{star})
\begin{eqnarray*}
\left\{
\begin{array}{l}
{\mathcal L}(X,D_X)\,{\mathcal U}={\mathcal F}
\qquad\mbox{in}\,\,\Omega,
\\[5pt] 
{\displaystyle{\frac{\partial^k{\mathcal U}}{\partial\nu^k}}}
=g_k\quad\,\,\mbox{on}\,\,\partial\Omega,\,\,\,\,\,0\leq k\leq m-1.
\end{array}
\right.
\end{eqnarray*}

\medskip

\begin{proposition}\label{pr2}
 If ${\mathcal U}\in W_p^{m,a}(\Omega)$ then, 
necessarily, 
$${\mathcal F}\in V_p^{-m,a}(\Omega), \,\,\, 
g:=\{g_k\}_{0\leq k\leq m-1}\in \dot{W}^{m-1+s}_p(\partial\Omega)$$
and, moreover, 
\begin{eqnarray*}
\|g\|_{\dot{W}^{m-1+s}_p(\partial\Omega)}
+\|{\mathcal F}\|_{V_p^{-m,a}(\Omega)}
\leq C\|{\mathcal U}\|_{W_p^{m,a}(\Omega)}.
\end{eqnarray*}
\end{proposition}

For the results in  the converse direction,
the main hypothesis in \cite{[MMS1]} 
requires that, {\it at small scales, the 
so called local  mean oscillations of the unit normal to $\partial\Omega$ and of 
the coefficients of the differential operator ${\mathcal L}(X,D_X)$ are 
not too large, relative to the Lipschitz constant of the domain $\Omega$,
the ellipticity constant of ${\mathcal L}(X,D_X)$, and the indices of the 
corresponding Besov space}.

\medskip

By the  {local mean oscillation}  of 
$F\in L_1(\Omega)$ we  understand 
\begin{equation*}
 \{F\}_{{\rm Osc}(\Omega)}:=
 \mathop{\hbox{lim}}_{\varepsilon\to 0} 
\Bigl(\mathop{\hbox{sup}}_{{\{B_\varepsilon\}}_\Omega}
\meanint_{\!\!\!B_\varepsilon\cap\Omega}\,\,
\meanint_{\!\!\!B_\varepsilon\cap\Omega}\,
\Bigl|\,F(x)-F(y)\,\Bigr| dxdy \Bigr), 
\end{equation*}
\noindent  where  $\{B_\varepsilon\}_\Omega$  stands for the  family of 
balls of radius $\varepsilon$ centered at points of $\Omega$ and the barred integrals denote the mean values. 
 Similarly, the local mean oscillation of $f\in L_1(\partial\Omega)$ is 
\begin{equation*}
 \{f\}_{{\rm Osc}(\partial\Omega)}:=
 \mathop{\hbox{lim}}_{\varepsilon\to 0}
\Bigl(\mathop{\hbox{sup}}_{\{B_\varepsilon\}_{\partial\Omega}}
\meanint_{\!\!\!B_\varepsilon\cap\partial\Omega}\,\,
\meanint_{\!\!\!B_\varepsilon\cap\partial\Omega}\,
\Bigl|\,f(x)-f(y)\,\Bigr|ds_xds_y \Bigr), 
\end{equation*}
\noindent  where 
$\{B_\varepsilon\}_{\partial\Omega}$  is the collection 
of $n$-dimensional balls of radius 
$\varepsilon$ with centers on $\partial\Omega$. 

\smallskip

Note that smallness of the local mean oscillation $\{\nu\}_{{\rm Osc}(\partial\Omega)}$ does not imply smallness of the Lipschitz constant. Indeed, let
$$\Omega = \{(x,y)\in \Bbb{R}^2, \, y>\varphi_\varepsilon (x)\},$$
where
$$\varphi_\varepsilon (x) = x\, \sin (\varepsilon \log |x|^{-1}).$$
Then $\|\varphi'_\varepsilon\|_{L_\infty(\Bbb{R})} \sim 1$, while  $\|\varphi'_\varepsilon\|_{BMO(\Bbb{R})} \leq C\, \varepsilon$.

\subsection{Solvability of the Dirichlet problem in $W_p^{m,a}$}
 
The main result in \cite{[MMS1]} runs as follows.

  Let $\Omega\subset\Bbb{R}^n$ be a bounded Lipschitz domain whose Lipschitz 
constant is does not exceed $M$, and assume that the operator ${\mathcal L}(X,D_X)$ of order $2m$
is strongly elliptic, and has bounded, measurable complex  coefficients.
 
\begin{theorem}\label{th1}
There exists  a positive constant $C$, depending only on $M$ and the 
ellipticity constant of ${\mathcal L}$, such that:
For each $p\in (1,\infty)$, $s\in (0,1)$ and $a:=1-s-1/p$, the Dirichlet 
problem 
\begin{eqnarray}\label{eq1a}
\left\{
\begin{array}{l}
\displaystyle{\sum_{|\alpha|=|\beta|=m}
D^\alpha(A_{\alpha\beta}(X)\,D^\beta\,{\mathcal U})}={\mathcal F} 
\qquad\mbox{for}\,\,X\in\Omega,
\\[20pt] 
{\displaystyle\frac{\partial^k{\mathcal U}}{\partial\nu^k}}=g_k
\,\,\quad\mbox{on}\,\,\partial\Omega,\qquad 0\leq k\leq m-1.
\end{array}
\right.
\end{eqnarray}
with ${\mathcal F} \in V_p^{-m,a}(\Omega)$ and $g:=\{g_k\}_{0\leq k\leq m-1}$ in $\dot{W}^{m-1+s}_p(\partial\Omega)$
has a unique solution 
${\mathcal U}\in W^{m,a}_p(\Omega)$ if 
 the coefficient matrices $A_{\alpha\beta}$ and the 
exterior normal vector $\nu$ to $\partial\Omega$ satisfy 
\begin{eqnarray}\label{a0}
&& \{\nu\}_{{\rm Osc}(\partial\Omega)}
+\sum_{|\alpha|=|\beta|=m}\{ A_{\alpha\beta}\}_{{\rm Osc}(\Omega)}
\\
&& \leq\,C\,s(1-s)\Bigl(p^2(p-1)^{-1}+s^{-1}(1-s)^{-1}\Bigr)^{-1}.\nonumber
\end{eqnarray}
For second order operators the factor $s(1-s)$ 
 in the last inequality 
can be removed.
Furthermore, there exists 
$C=C(\partial\Omega,{\mathcal A},p,s)>0$ 
such that 
\begin{equation}\label{estUU}
\|{\mathcal U}\|_{W_p^{m,a}(\Omega)}
\leq C\Big(\|g\|_{\dot{W}^{m-1+s}_p(\partial\Omega)}
+\|{\mathcal F}\|_{V_p^{-m,a}(\Omega)}\Bigr). 
\end{equation}

\end{theorem}

\medskip

The next assertion, obtained in \cite{[MMS1]}, is a byproduct of the proof of  Theorem \ref{th1}.

\begin{theorem}\label{th22}
 Let $\Omega\subset\Bbb{R}^n$ be a bounded Lipschitz domain whose Lipschitz 
constant does not exeed $ M$, and assume that the operator ${\mathcal L}(X,D_X)$ of order $2m$
is strongly elliptic, and has bounded, measurable (complex) coefficients.
 Then there exists $\varepsilon$, depending only on $M$ and the 
ellipticity constant of ${\mathcal L}$, such that:

For each $p\in (1,\infty)$, $s\in (0,1)$ and $a:=1-s-1/p$, the Dirichlet 
problem for ${\mathcal L}$  with ${\mathcal F} \in V_p^{-m,a}(\Omega)$ and $g:=\{g_k\}_{0\leq k\leq m-1}$ in $\dot{W}^{m-1+s}_p(\partial\Omega)$
has a unique solution 
${\mathcal U}\in W^{m,a}_p(\Omega)$ if 
\begin{eqnarray*}
|2^{-1}-p^{-1}|<\varepsilon\quad\mbox{ and }\quad|a|<\varepsilon,
\end{eqnarray*}
\end{theorem}
 
Recently M. Agranovich \cite{[Ag1]} obtained  this type of results ($|p-2|$ is small and $0<s<1$)  for  
 both Dirichlet and Neumann problems for a subclass of strongly elliptic systems with Douglis-Nirenberg structure in bounded Lipschitz domains. The regularity results in \cite{[Ag1]} concern solutions in spaces of Bessel potentials $H_p^\sigma$ and Besov spaces $B_p^\sigma$ with coefficients of differential operators satisfying the uniform Lipschitz condition. The approach is based on regularity methods due to Savar\'e as well as on author's developement of interpolation theory of spaces $H_p^\sigma$ and $B_p^\sigma$ with $\sigma$ of arbitrary sign (see \cite{[Ag2]} for more details), where essential role is played by an extension operator from $\Omega$ to ${\Bbb{R}}^n$ introduced by V.\,Rychkov \cite{[Ry]}. To be specific, Agranovich considered a Douglis-Nirenberg system 
with the principal part  $L_0$ whose entries are given by
$$L_{j,k} (x, D) = \sum_{|\alpha| =m_j,\atop{|\beta| =m_k}} D^\alpha\bigl( a_{\alpha,\beta}^{j,k}(x) \, D^\beta \bigr)$$
The coefficients $a_{\alpha,\beta}^{j,k}(x)$ are complex-valued and the formal self-adjointness of the operatot $L_{j,k}$ is not assumed. The principal symbol of the system, i. e. the matrix $L_0(x,\xi)$ with entries $L_{j,k}(x,\xi)$ is subordinate to the condition of strong ellipticity:
$${\rm Re}\, L_0 (x, \xi) \geq C\, \Lambda(\xi),$$
where $\Lambda(\xi)$ is the diagonal matrix with entries $|\xi|^{2m_j} $ on the main diagonal, $ j= 1, \ldots, l$, and $C$ is a positive constant. The use the Savar\'e method requires the additional condition 
$${\rm Re} \sum_{j,k =1}^l\sum_{|\alpha| =m_j,\atop{|\beta| =m_k}}a_{\alpha,\beta}^{j,k}(x) \zeta_\beta^k\, \overline{\zeta_\alpha^j}\geq 0$$
for any numbers $\zeta_{\alpha}^k$ at all points $x\in \Omega$.

\medskip

Let  ${\rm BMO}$  and  ${\rm VMO}$ 
stand, respectively, for the   space of functions 
of bounded mean oscillations and the  space 
of functions of vanishing mean oscillations (considered either on $\Omega$, 
or on $\partial\Omega$). It can be proved that 
$$
\{F\}_{{\rm Osc}}\sim{\rm dist}\,(F,{\rm VMO})
$$ 
where the distance is taken in ${\rm BMO}$. Thus
the small oscillation condition introduced in  Theorem \ref{th1}  holds if 
\begin{eqnarray*}
&& {\rm dist}\,(\nu,{\rm VMO})
+\sum_{|\alpha|=|\beta|=m}{\rm dist}\,(A_{\alpha\beta},{\rm VMO})
\\[4pt]
&& \leq\,C\,s(1-s)\Bigl(p^2(p-1)^{-1}+s^{-1}(1-s)^{-1}\Bigr)^{-1}. 
\end{eqnarray*}
\noindent  This is the case if, e.g.,  
 $[\nu]_{\rm BMO}+\sum[A_{\alpha\beta}]_{\rm BMO}$ 
is sufficiently small hence, trivially, if 
$\nu\in {\rm VMO}(\partial\Omega)$ and  $A_{\alpha\beta}$ belong 
to ${\rm VMO}(\Omega)$, irrespective of $p$, $s$, 
${\mathcal L}$ and $\Omega$. 

\smallskip

 Other examples of domains satisfying the hypotheses of Theorem \ref{th1} are: 
 Lipschitz domains with a sufficiently small Lipschitz constant, 
relatively to the exponents $p$ and $s$.
 In particular, 
Lipschitz polyhedral domains with dihedral angles sufficiently close, 
depending on $p$ and $s$,  to $\pi$.

\medskip

 The innovation in \cite{[MMS2]} that allows  to consider 
  boundary data in higher-order Besov spaces, is the 
systematic use of   weighted Sobolev spaces. 
In relation to the standard Besov scale, 
we would like to point out that 
\begin{equation*}
a=1-s-\frac{1}{p}\in (0,1-1/p)\Longrightarrow 
W^{m,a}_p(\Omega)\hookrightarrow B^{m-1+s+1/p}_{p}(\Omega),
\end{equation*}
{ and } 
\begin{equation*}
a=1-s-\frac{1}{p}\in (-1/p,0)\Longrightarrow
 B^{m-1+s+1/p}_{p}(\Omega)\hookrightarrow W^{m,a}_p(\Omega).
\end{equation*}
\noindent  Of course, $W^{m,a}_p(\Omega)$ is just the classical Sobolev 
space $W^{m}_p(\Omega)$ when $a=0$. 

\medskip

  In the context of   Theorem \ref{th1}, for  a Lipschitz $\partial\Omega$,  
\begin{equation*}
 \sum_{|\alpha|\leq m-1}
\|{\rm Tr}\,[D^\alpha\,{\mathcal U}]\|_{B_{p}^{s}(\partial\Omega)}
 \sim \Bigl(\sum_{|\alpha|\leq m}\int_{\Omega}\rho(X)^{p(1-s)-1}\,
|D^\alpha{\mathcal U}(X)|^p\,dX\Bigr)^{1/p},
\end{equation*}
\noindent  uniformly in ${\mathcal U}$ satisfying 
 ${\mathcal L}(X,D_X)\,{\mathcal U}=0$ in $\Omega$. 

\smallskip

This generalizes the trace and extension  result  by S. Uspenski\u{\i} \cite{[Usp]}
 according to which 
$$
\|{\rm Tr}\,u\|_{B_{p}^s(\Bbb{R}^{n-1})}\sim
\Bigl(\int_{\Bbb{R}^n_+}x_n^{p(1-s)-1}|\nabla u(x',x_n)|^p
dx\Bigr)^{1/p}
$$
 if $1<p<\infty$ and $0<s<1$, uniformly for $u$  harmonic 
in the upper-half space. 

\smallskip

 Of course, condition (\ref{A-bdd}) ensures that the left-hand side of 
(\ref{a0}) is always finite but it is its actual size which determines
whether for a given pair of indices $s$, $p$, the problem 
(\ref{eq1}), (\ref{1.6}), (\ref{1.9}) is well-posed. Note that 
the maximum value that the right-hand side of (\ref{a0}) takes for
$0<s<1$ and $1<p<\infty$ occurs precisely when $p=2$ and $a:=1-s-1/p=0$. 
As (\ref{a0}) shows, the set of pairs $(s,1/p)\in(0,1)\times(0,1)$ for which 
(\ref{eq1}) is well-posed in the context of Theorem~\ref{th1} 
exhausts the entire square $(0,1)\times(0,1)$ as the distance from 
$\nu$ and the $A_{\alpha\beta}$'s to {\rm VMO} tends to zero 
(while the Lipschitz constant of $\Omega$ and the ellipticity constant of
${\mathcal L}$ stay bounded). That the geometry of the Lipschitz domain 
$\Omega$ intervenes in this process through a condition such as (\ref{a0}) 
confirms a conjecture made 
by P.\,Auscher and M.\,Qafsaoui in {{\cite{[AQ]}}}.

\smallskip

Theorem \ref{th22} can be viewed as a far reaching  extension of a 
well-known theorem of N.\,Meyers, who has treated the case $m=1$, $l=1$ 
in {{\cite{[Mey]}}}. The example given in Section 5 of {{\cite{[Mey]}}} 
shows that the inclusion  of $p$ into a small neighborhood of $2$ is a 
necessary condition, even when $\partial\Omega$ is smooth, if 
the coefficients $A_{\alpha\beta}$ are merely bounded. For higher 
order operators we make use of an example  due to V.\,Maz'ya 
{{\cite{[Maz1]}}} (cf. also the contemporary article by E.\,De Giorgi 
{{\cite{[DG2]}}}). Specifically, when 
$m\in\NN$ is even, consider the divergence-form equation
\begin{equation}\label{Maz-Op}
\Delta^{\frac{1}{2}m-1}{\mathcal L}_4\,\Delta^{\frac{1}{2}m-1}{\mathcal U}=0
\quad\mbox{in }\,\,\Omega:=\{X\in\RR^n:\,|X|<1\},
\end{equation}
\noindent where ${\mathcal L}_4$ is the fourth order operator  
\begin{eqnarray}\label{Maz-Op2}
{\mathcal L}_4(X,D_X)\,{\mathcal U} & := & 
a\,\Delta^2{\mathcal U}
+b\sum_{i,j=1}^n\Delta\Bigl(\frac{X_iX_j}{|X|^2}
\partial_i\partial_j\,{\mathcal U}\Bigr)
+b\sum_{i,j=1}^n\partial_i\partial_j\Bigl(\frac{X_iX_j}{|X|^2}\,
\Delta\,{\mathcal U}\Bigr)
\nonumber\\[6pt]
&&+c\sum_{i,j,k,l=1}^n\partial_k\partial_l\Bigl(\frac{X_iX_jX_kX_l}{|X|^4}
\partial_i\partial_j
\,{\mathcal U}\Bigr).
\end{eqnarray}
\noindent Obviously, the coefficients of ${\mathcal L}_4(X,D_X)$ are 
bounded, and if the parameters $a,b,c\in\RR$, $a>0$, are chosen such 
that $b^2<ac$ then ${\mathcal L}$ along with 
$\Delta^{\frac{1}{2}m-1}{\mathcal L}_4\,\Delta^{\frac{1}{2}m-1}$ 
are strongly elliptic. Now, if $W^s_p$ denotes the usual $L_p$-based 
Sobolev space of order $s$, it has been observed in {{\cite{[Maz1]}}} 
that the function ${\mathcal U}(X):=|X|^{\theta+m-2}\in W^m_2(\Omega)$ has 
${\rm Tr}\,{\mathcal U}\in C^\infty(\partial\Omega)$ and is a
weak solution of (\ref{Maz-Op}) for the choice 
\begin{equation}\label{Maz-Op3}
\theta:=2-\frac{n}{2}+\sqrt{\frac{n^2}{4}-\frac{(n-1)(bn+c)}{a+2b+c}}.
\end{equation}
\noindent Thus, if 
$$a:=(n-2)^2+\varepsilon, \,\,\,  b:=n(n-2), \,\,\,  c:=n^2, \,\,\, 
\varepsilon>0,$$
 the strong ellipticity condition is satisfied and 
$\theta=\theta(\varepsilon)$ becomes 
$$2-n/2+n\,\varepsilon^{1/2}/2\bigl(4(n-1)^2+\varepsilon\bigr)^{1/2}.$$ 
However, ${\mathcal U}\in W^m_p(\Omega)$ if and only if 
$p<n/(2-\theta(\varepsilon))$, and the bound $n/(2-\theta(\varepsilon))$ 
approaches $2$ when $\varepsilon\to 0$. An analogous example can be produced 
when $m>1$ is odd, starting with a sixth order operator 
${\mathcal L}_6(X,D_X)$ from {{\cite{[Maz1]}}}. 
In the above context, given that 
$W^1_n(\Omega)\hookrightarrow{\rm VMO}(\Omega)$, it is significant to point 
out that both for the example in {{\cite{[Mey]}}}, when $n=2$, and for 
(\ref{Maz-Op}) when $n\geq 3$, the coefficients have their 
gradients in weak-$L_n$,  yet they fail to belong to $W^1_n(\Omega)$.

\smallskip

We conclude this subsection with a remark pertaining to the presence of lower 
order terms. More specifically, granted Theorem~\ref{th1}, a standard 
perturbation argument (cf., e.g., \cite{[H]}) proves the following. 
Assume that 
\begin{equation}\label{E444-bis}
{\mathcal A}(X,D_X)\,{\mathcal U}
:=\sum_{0\leq |\alpha|,|\beta|\leq m}D^\alpha({\mathcal A}_{\alpha\beta}(X)
\,D^\beta{\mathcal U}),\qquad X\in\Omega,
\end{equation}

\noindent 
where the principal  part of ${\mathcal A}(X,D_X)$ satisfies the 
hypotheses made in Theorem~\ref{th1} and the lower order terms are bounded.
Then, assuming that  (\ref{a0}) holds, the 
Dirichlet problem (\ref{eq1a}) is Fredholm with index zero, in the sense that 
the operator 
$$W^{m,a}_p(\Omega)\ni{\mathcal U}\mapsto
\Bigl({\mathcal A}(X,D_X)\,{\mathcal U}\,\,,\,\,
\{\partial^k{\mathcal U}/\partial\nu^k\}_{0\leq k\leq m-1}\Bigr)
\in V_p^{-m,a}(\Omega)\oplus\dot{W}^{m-1+s}_p(\partial\Omega)$$
\noindent is so. Furthermore, the estimate 
\begin{equation}\label{estUU-bis}
\|{\mathcal U}\|_{W_p^{m,a}(\Omega)}
\leq C\, \Bigl(\|{\mathcal F}\|_{V_p^{-m,a}(\Omega)} 
+\|g\|_{\dot{W}^{m-1+s}_p(\partial\Omega)}+\|{\mathcal U}\|_{L_p(\Omega)}\Bigr)
\end{equation}
\noindent holds for any solution ${\mathcal U}\in W_p^{m,a}(\Omega)$
of (\ref{eq1a}).

\subsection{Comments on the proof of  Theorem \ref{th1}}

 One difficulty linked with the case $m>1$ arises from the 
way the norm in 
\begin{equation}\label{W-Nr}
\Bigl(\sum_{|\alpha|\leq m}\int_\Omega|D^\alpha{\mathcal U}(X)|^p 
\rho(X)^{pa}\,dX\Bigr)^{1/p}<\infty.
\end{equation}
 behaves under a change of variables
$$\varkappa:\Omega=\{(X',X_n):\,X_n>\varphi(X')\}\to\RR^n_+$$
 destined to 
flatten the Lipschitz surface $\partial\Omega$. When $m=1$, a simple 
bi-Lipschitz changes of variables such as 
$\Omega\ni (X',X_n)\mapsto(X',X_n-\varphi(X'))\in\RR^n_+$ will do, but 
matters are considerable more subtle in the case $m>1$. The 
 extension operator  used in \cite{[MMS1]} was introduced by J.\,Ne\v{c}as (in a different context; cf. p.\,188 in {\cite{[Nec]}}) 
and then  rediscovered  by V.\,Maz'ya and 
T.\,Shaposhnikova in {\cite{[MS1]}} (see also \cite{[MS2]}, and later by  B.\,Dahlberg, C.\,Kenig J.\,Pipher, E.\,Stein and G.\,Verchota 
(cf. {\cite{[Dah2]}} and the discussion in {\cite{[DKPV]}}), 
and S.\,Hofmann and J.\,Lewis in {\cite{[HL]}}. 

\smallskip

The extension operator in question is introduced in the following way. 
Fix a smooth, radial, 
decreasing, even, non-negative function $\zeta$ in $\mathbb{R}^{n-1}$ 
such that $\zeta(t)=0$ for $|t|\geq 1$ and 
\begin{equation}\label{zeta-int}
\int\limits_{\mathbb{R}^{n-1}}\zeta(t)\,dt=1.
\end{equation}
\noindent (For example,  
$\zeta(t):=c\,{\rm exp}\,(-1/(1-|t|^2)_+)$ for a suitable $c$.) 
Define the extension operator $T$ by 
\begin{equation}\label{10.1.20}
(T\varphi)(x',x_n):=\int\limits_{\mathbb{R}^{n-1}}\zeta(t)\varphi(x'+x_nt)\,dt,
\qquad(x',x_n)\in\RR^n_+,
\end{equation}
\noindent acting on functions $\varphi$ from $L_{1 ,loc}(\mathbb{R}^{n-1})$. 

\smallskip

The two estimates below provide  useful properties of the  operator $T$.

{\rm (i)} For each multi-indices $\alpha$ with $|\alpha|>1$
there exists $c>0$ such that 
\begin{equation*}
\Bigl|D^\alpha_{x}(T\varphi)(x)\Bigr|
\leq c\,x_n ^{1-|\alpha|}[\nabla\varphi]_{\rm BMO(\mathbb{R}^{n-1})},
\quad\forall\,x=(x',x_n)\in\mathbb{R}^n_+. 
\end{equation*}

{\rm (ii)} If $\nabla_{x'}\varphi\in{\rm BMO}(\mathbb{R}^{n-1})$ then 
$\nabla(T\varphi)\in{\rm BMO}(\mathbb{R}^n_+)$ and 
\begin{equation*}
[\nabla(T\varphi)]_{{\rm BMO}(\mathbb{R}^n_+)}
\leq c\,[\nabla_{x'}\varphi]_{{\rm BMO}(\mathbb{R}^{n-1})}.
\end{equation*}

\medskip

Another ingredient  in the proof of Theorem \ref{th1} is deriving estimates for 
$D_x^\alpha D_y^\beta G(x,y)$ where $G$ is the Green matrix of the operator 
$${L}(D_x) = \sum_{|\alpha|=2m}A_\alpha D_x^\alpha,$$ 
 i.e. a unique solution of the boundary-value problem
$$\left\{
\begin{array}{l}
{L}(D_x)G(x,y)=\delta(x-y)I_l\quad\mbox{for}\,\,x\in\mathbb{R}^n_+,
\\[2pt]
\displaystyle{\Bigl(\frac{\partial ^j}{\partial x_n^j}G\Bigr)((x',0),y)}
=0\,I_l\qquad
\mbox{for}\,\,x'\in\mathbb{R}^{n-1}, \,\,\, 0\leq j\leq m-1,
\end{array}
\right.$$
\noindent where $y\in\mathbb{R}^n_+$ is regarded as a parameter. The methods employed in earlier 
works are  based on explicit representation formulas for $G(x,y)$ and 
 cannot be adapted  to the case of  non-symmetric, complex 
coefficient, higher order systems. The  approach in \cite{[MMS1]} 
consists of proving  that the residual part 
$R(x,y):=G(x,y)-\Phi(x-y)$, where $\Phi$ is a fundamental matrix  
for $L(D_x)$, has the property 
\begin{equation*}
\|D^\alpha_xD^\beta_y R(x,y)\|
\leq C\,|x-\bar{y}|^{-n}
\end{equation*}

\noindent 
for $|\alpha|=|\beta|=m$, $x,y\in\RR^n_+$, where $\bar{y}:=(y',-y_n)$ is
the reflection of the point $y\in\RR^n_+$ with respect to $\partial\RR^n_+$.

\section{Stokes system}

By $B^s_{p, q}(\mathbb{R}^n)$ we denote the space of functions in $\mathbb{R}^n$ having the finite norm 
\begin{equation}\label{eq13.8u}
 \| u\|_{B^s_{p, q}(\mathbb{R}^n)} =\Bigl( \int_{\mathbb{R}^n} \| \Delta_h \nabla_{[s]} u \|^q_{L_p(\mathbb{R}^n)} |h|^{-n-q\{s\}}\,d h \Bigr)^{1/q}+\|u\|_{W^{[s]}_p(\mathbb{R}^n)}\;,
 \end{equation}
where $\{s\}>0$, $p$, $q \ge 1$ and $\Delta_h v = v(\cdot + h) - v(\cdot)$. For $p=q$ we use the notation $B_p^s(\mathbb{R}^n)$.

\smallskip

The Besov  scale $B^s_{p, q}(\Omega)$ 
is defined by restricting the (tempered) distributions from the corresponding
spaces in ${\mathbb{R}}^n$ to the open set $\Omega$. 
Also, $B_{p,q}^s(\partial\Omega)$ stands for the Besov class on the Lipschitz 
manifold $\partial\Omega$, obtained by transporting (via a partition of unity 
and pull-back) the standard scale $B_{p,q}^s({\mathbb{R}}^{n-1})$.

\subsection{Weak solvability in Besov and Triebel-Lizorkin spaces}
 
Consider the Stokes system in an arbitrary bounded Lipschitz domain $\Omega\subset{\mathbb{R}}^n$, 
$n\geq 2$,

\begin{eqnarray}\label{eq1.24}
\begin{array}{l}
\Delta {u}-\nabla\pi= {f}\in B^{s+\frac{1}{p}-2}_{p,q}(\Omega),
\quad{\rm div}\, {u}=g\in B^{s+\frac{1}{p}-1}_{p,q}(\Omega),
\\[10pt]
 {u}\in B^{s+\frac{1}{p}}_{p,q}(\Omega),\quad
\pi\in B^{s+\frac{1}{p}-1}_{p,q}(\Omega),\quad
{\rm Tr}\, {u}= {h}\in B^s_{p,q}(\partial\Omega),
\end{array}
\end{eqnarray}

subject to the (necessary) compatibility condition
\begin{eqnarray}\label{comp-CD}
\int_{\partial{\mathcal{O}}}\langle\nu, {h}\rangle\,d\sigma
=\int_{\mathcal{O}}g(X)\,dX, 
\quad\mbox{for every component ${\mathcal{O}}$ of $\Omega$}. 
\end{eqnarray}

\medskip

 When $\partial\Omega$ is 
sufficiently smooth (at least of the class $C^2$), the problem (\ref{eq1.24}) was studied in many papers, first in  Sobolev spaces 
with an integer amount of smoothness by 
V.A.\,Solonnikov \cite{[Sol]}, L.\,Cattabriga \cite{[Cat]}, R.\,Temam \cite{[Tem]}, 
Y.\,Giga \cite{[Gig]},  
R.\,Dautray and J.-L.\,Lions \cite{[DL]}, among others. This has been subsequently extended by 
C.\,Amrouche and V.\,Girault \cite{[AG]} to the case when 
$\partial\Omega\in C^{1,1}$ and, further, by G.P.\,Galdi, C.G.\,Simader 
and H.\,Sohr \cite{[GSS]} when $\partial\Omega$ is Lipschitz, with a small  
Lipschitz constant. The case when $\Omega$ is a polygonal domain in 
${\mathbb{R}}^2$, or a polyhedral domain in ${\mathbb{R}}^3$ also has a reach history. An extended 
account of this field of research can be found in V.\,Kozlov, V.\,Maz'ya 
and J.\,Rossmann's monograph \cite{[KMR]}, which also contains 
references to earlier work. Among recent publications we mention the paper 
by V.\,Maz'ya and J.\,Rossmann \cite{[MR1]} as well as their book \cite{[MR2]}. Lipschitz and $C^1$ subdomains of Riemannian 
manifolds were treated in  the paper \cite{[DMi]} by M.\,Dindo\v{s} and M.\, Mitrea, 
and the paper \cite{[MT]} by M.\, Mitrea and M.\,Taylor. 

\smallskip

The principal result in \cite{[MMS2]} on solutions of the problem (\ref{eq1.24}) holds under the mild condition on the normal $\nu$ to the boundary of the Lipschitz graph domain and 
runs as follows.

\begin{theorem}\label{Th2}
Assume that 
\begin{equation}\label{IND}
\frac{n-1}{n}<p\leq\infty, \qquad 0<q\leq\infty, \qquad (n-1)\bigl({\textstyle{\frac{1}{p}}}-1\bigr)_+<s<1.
\end{equation}
 Then there exists $\delta>0$ which depends only on the 
Lipschitz character of $\Omega$ and the exponent $p$, with the property
that if $\{\nu\}_{{\rm Osc}(\partial\Omega)} <\delta$, then the problem $(\ref{eq1.24})$ is well-posed (with uniqueness modulo locally constant 
functions in $\Omega$ for the pressure). There exists a finite, positive
constant $C=C(\Omega,p,q,s,n)$ such that
$$
\| {u}\|_{B^{s+\frac{1}{p}}_{p,q}(\Omega)}
+\inf_{c}\|\pi-c\|_{B^{s+\frac{1}{p}-1}_{p,q}(\Omega)}$$ 
$$\leq C\Bigl(\| {f}\|_{B^{s+\frac{1}{p}-2}_{p,q}(\Omega)}
+ \|g\|_{B^{s+\frac{1}{p}-1}_{p,q}(\Omega)}
+ \| {h}\|_{B^s_{p,q}(\partial\Omega)}\Bigr),
$$
where the infimum is taken over all locally constant functions 
$c$ in $\Omega$. 
\end{theorem}

Moreover, analogous well-posedness result holds on the Triebel-Lizorkin scale (for its definition see, for example, the book by T. Runst and W. Sickel \cite{[RS]}), 
i.e. for the problem
\begin{equation}\label{eq1.26}
\begin{array}{l}
\Delta {u}-\nabla\pi= {f}\in F^{s+\frac{1}{p}-2}_{p,q}(\Omega),\quad
{\rm div}\, {u}=g\in F^{s+\frac{1}{p}-1}_{p,q}(\Omega),
\\[10pt]
 {u}\in F^{s+\frac{1}{p}}_{p,q}(\Omega),\quad
\pi\in F^{s+\frac{1}{p}-1}_{p,q}(\Omega),\quad
{\rm Tr}\, {u}= {h}\in B^s_{p,q}(\partial\Omega),
\end{array}
\end{equation}
This time, in addition to the previous conditions imposed on the indices
$p$, $q$, it is also assumed that $p,q<\infty$.

\medskip

It should be noted that 
conditions (\ref{IND}) describe the largest range of indices $p,q,s$ for 
which the Besov spaces $B^s_{p,q}(\partial\Omega)$ can be meaningfully 
defined on the Lipschitz manifold $\partial\Omega$.

\smallskip

Regarding the nature of the main result, note that
 no topological restrictions are imposed on the Lipschitz 
domain $\Omega$ (in particular, the boundary can be disconnected). 
This is significant since the approach is via boundary layer
potentials, whose invertibility properties are directly affected by 
topological nature of the domain.  That difficulty is overcome by  using certain  mapping properties of the hydrostatic layer potentials: single layer potential ${\mathcal{S}}$ and double layer potential ${\mathcal{D}}_{\lambda}$  for the velocity and single layer potential ${\mathcal{Q}}$ and double layer potential ${\mathcal{P}}_{\lambda}$  for the presure. We list their properties in the following theorem.

\begin{theorem}\label{P-Lay}
Let $\Omega$ be a bounded Lipschitz domain in ${\mathbb{R}}^n$, $n\geq 2$, 
and assume that $\lambda\in{\mathbb{R}}$, $\frac{n-1}{n}<p\leq\infty$, 
$(n-1)(\frac{1}{p}-1)_+<s<1$, and $0<q\leq\infty$. Then
\begin{eqnarray*}
&& {\mathcal{S}}:B^{s-1}_{p,q}(\partial\Omega)
\longrightarrow B^{s+\frac{1}{p}}_{p,q}(\Omega),
\\[4pt]
&& {\mathcal{D}}_{\lambda}:B^s_{p,q}(\partial\Omega)
\longrightarrow B^{s+\frac{1}{p}}_{p,q}(\Omega),
\\[4pt]
&& {\mathcal{Q}}:B^{s-1}_{p,q}(\partial\Omega)
\longrightarrow B^{s+\frac{1}{p}-1}_{p,q}(\Omega),
\\[4pt]
&& {\mathcal{P}}_{\lambda}:B^s_{p,q}(\partial\Omega)
\longrightarrow B^{s+\frac{1}{p}-1}_{p,q}(\Omega),
\end{eqnarray*}
\noindent are well-defined, bounded operators. 
\end{theorem}

\subsection{Dirichlet data in $L_p(\Omega)$}
We  say a few words about the Dirichlet problem for the multi-dimensional Stokes system with $L_p$ boundary data which is not touched upon in \cite{[MMS1]}. Following the paper by M. Dindo\v{s}  and V. Maz'ya \cite{[DMa]}, we shall speak  about both the Lam\'e system (with the Poisson ratio $\alpha <1/2$) and the Stokes system (with $\alpha =1/2$).

\smallskip

Let us consider a bounded domain $\Omega$ in ${\mathbb{R}}^n$, $n\geq 3$, and the system
\begin{equation}\label{n1}
\Delta {u}-\nabla\pi= 0,\quad
{\rm div}\, {u} + (1-2\alpha)\pi =0 \quad {\rm in} \,\, \Omega
\end{equation}
complemented by the condition
\begin{equation}\label{n2}
{\rm Tr}\, {u}=  h \in L_p(\partial\Omega),
\end{equation}
and the class of solutions is described by the inclusion
\begin{equation}\label{n3}
u^* \in L_p(\partial\Omega).
\end{equation}
The boundary values ${\rm Tr}\, {u}$ are understood in nontangential sense, that is in the sense of the limit
$${\rm Tr}\, {u} (x) = \lim_{y\to x, y\in \Gamma(x)} u(y),$$
over a collection of interior nontangential cones $\Gamma(x)$ of some aperture and height and vertex at $x\in \partial\omega$, and $u^*$ is the classical nontangential maximal function defined as
$$u^*(x) = \sup_{y\in \Gamma(x)} |u(y)| \qquad {\rm for}\,\, {\rm all} \,\, x\in \partial\omega.$$

\smallskip

One says that  problem (\ref{n1}), (\ref{n2}) is $L_p$ solvable if for all vector fields $h\in L_p(\partial\Omega)$ there is a pair $(u,\pi)$ satisfying (\ref{n1}) - (\ref{n3}) and moreover, for some $C>0$ independent of $h$, the estimate
$$\|u^*\|_{L_p(\partial\Omega)} \leq C\, \|h\|_{L_p(\partial\Omega)}$$
holds. Furthermore, the problem (\ref{n1}), (\ref{n2}) is said to be solvable for continuous data if, for all $h\in C(\partial\Omega)$ the vector field $u$ belongs to $C(\overline\Omega)$ and the estimate
$$\|u^*\|_{C(\overline\Omega)} \leq C\, \|h\|_{C(\partial\Omega)}$$
holds. 

\smallskip

The $L_p$ solvability of the Dirichlet problem for the Stokes system is established by Z. Shen \cite{[Sh]} 
for all $p\in (2-\varepsilon(\Omega), \infty]$ provided $\Omega$ is a three-dimensional Lipschitz graph domain. In \cite{[DMa]}, the problem (\ref{n1}), (\ref{n2}) is considered on domains in ${\mathbb{R}}^n$, $n\geq 3$, with isolated conical singularity (not necessary a Lipschitz graph) and the authors prove its solvability for all $p\in (2-\varepsilon(\Omega), \infty]$ as well as its solvability in $C(\overline\Omega)$ for the data in $C(\partial\Omega)$. This seems to be a strong indication that the range $p\in (2-\varepsilon(\Omega), \infty]$ should hold also for Lipschitz graph domains. However, there is no such a result for $n>3$.

\subsection{Lipschitz continuous solutions}

We cite a regularity result for solutions of the Dirichlet problem for the Stokes system in a plane convex domain obtained by V. Kozlov and V. Maz'ya in \cite{[KM3]}:
\begin{eqnarray}\label{n4}
&&-\Delta {u}+\nabla\pi= f,\quad
{\rm div}\, {u}=0 \quad {\rm in} \,\, \Omega\nonumber\\
\\
&& {\rm Tr}\, u=0 \qquad {\rm on} \,\,\, \partial\Omega,\nonumber
\end{eqnarray}
where $f\in W^{-1}_2(\Omega)$ and $(u, \pi) \in \ring W^{1}_2(\Omega) \times L_2(\Omega)$.

\begin{theorem}\label{T3}
Let $\Omega$ be a bounded convex two-dimensional domain and let $f\in L_q(\Omega)$ for some $q>2$. Then the  velocity vector $u\in \ring W_2^1(\Omega)$   admits the estimate
$$\|\nabla u\|_{L_\infty(\Omega)} \leq C \, \|f\|_{L_q(\Omega)},$$
where $C$ depends only on $\Omega$.
\end{theorem}

A direct consequence of this result for the nonlinear Navier-Stokes system is as follows.

\begin{corollary}\label{cor1}
Let $(u, \pi) \in \ring W^{1}_2(\Omega) \times L_2(\Omega)$ solve the Dirichlet problem
\begin{eqnarray}\label{n5}
&&-\Delta {u}+\nabla\pi +\sum_{k=1}^2 u_k \partial_{k} u 
= f,\quad
{\rm div}\, {u}=0 \quad {\rm in} \,\, \Omega,\nonumber\\
&& u=0 \qquad {\rm on} \,\,\, \partial\Omega,
\end{eqnarray}
where $f\in W^{-1}_2(\Omega)$. Let $\Omega$ be a bounded convex two-dimensional domain and let $f\in L_q(\Omega)$ for some $q>2$. Then the velocity vector $u\in \ring W_2^1(\Omega)$  belongs to the Lipschitz class $C^{0,1}(\overline\Omega)$.
\end{corollary}

\smallskip

A result of the same nature was obtained by V. Maz'ya in  \cite{[Maz3]} for solutions of the Neumann problem for the Poisson equation in arbitrary convex $n$-dimensional domain. Let $\Omega$ be a bounded convex  domain in ${\mathbb{R}}^n$ and let $W^{l}_p(\Omega)$ stand for the Sobolev space of functions in $L_p(\Omega)$ with distributional derivatives of order $l$ in $W^{l}_p(\Omega)$. By $L_{p,\perp}(\Omega)$ and $W^{l}_{p,\perp}(\Omega)$ one means the subpaces of functions in $L_p(\Omega)$ and $W^{l}_p(\Omega)$ subject to
$$\int_\Omega v\, dx =0.$$

\begin{theorem}\label{th5}
Let  $f\in L_{q,\perp}(\Omega)$ with a certain $q>n$ and let $u$ be the unique function in $W^{1}_2(\Omega)$, also orthogonal to $1$ in $L_2(\Omega)$ and satisfying the Neumann problem
\begin{eqnarray}\label{n6}
&&-\Delta {u}
= f \quad {\rm in}\,\, \Omega,\nonumber\\
\\
&& \frac{\partial u}{\partial \nu}=0 \quad {\rm on} \,\,\, \partial\Omega,\nonumber
\end{eqnarray}
where $\nu$ is the unit outward normal to $\partial\Omega$ and the problem $(\ref{n6})$ is understood in the variational sense. Then the solution $u\in W^{l}_{p,\perp}(\Omega)$ of the problem $(\ref{n6})$ satisfies the estimate
\begin{equation}\label{n7}
\|\nabla u\|_{L_\infty(\Omega)} \leq C\, \|f\|_{L_q(\Omega)}.
\end{equation}
\end{theorem}

As a particular case of a  result obtained by A. Cianchi and V. Maz'ya in \cite{[CM]} for a class of nonlinear equations, one can replace the $L_q$-norm on the right-hand side by the Lorentz norm $\|f\|_{L^{n,1}(\Omega)}$ which is the best possible majorant among those formulated in terms of Lorentz spaces. Similar facts for Lam\'e and Stokes  systems with boundary conditions different from those of Dirichlet are unknown.

\section{Higher regularity of solutions}

\subsection{Preliminaries}

In this section we give applications of Sobolev multipliers to the question of higher regularity in fractional Sobolev spaces of solutions to  boundary value problems for higher order elliptic equations  in   a Lipschitz domain.
 Since the sole Lipschitz graph property of $\Omega$ does not guarantee  higher regularity of solutions, we are forced to select an appropriate subclass of Lipschitz domains whose description involves a space of multipliers. For domains of this subclass we  develop a solvability and regularity theory analogous to the classical one for smooth domains.    We also show that the chosen  subclass of Lipschitz domains proves to be best possible in a certain sense.

\smallskip

Let $\Omega$ be a domain in $\Bbb{R}^n$ with compact closure. Throughout
this section we assume that for any point of the boundary
$\partial\Omega$ there exists a neighbourhood $U$ and a Lipschitz
function $\varphi$ such that 
\begin{equation}\label{eq1t}
U\cap \Omega = U \cap \{z = (x,y): x \in \Bbb{R}^{n-1},\; y >
\varphi (x)\}. 
\end{equation}

Let  $B_p^{l-1/p} (\Bbb{R}^{n-1}), l=1,2,\dots, $  denote 
the completion of the space $C^\infty_0 (\Bbb{R}^{n-1})$ in the norm
$$
\Bigl(\int_{\Bbb{R}^{n-1}}\int_{\Bbb{R}^{n-1}} \vert \nabla_{[l]-1} u (t) - \vert \nabla_{[l]-1} u (x)\vert^p
\vert x-t\vert^{-n+1-p\{l\}} dxdt\Bigr)^{1/p} + \Vert u \Vert_{L_p(\Bbb{R}^{n-1})}.
$$

Replacing  $\Bbb{R}^{n-1}$ by $\partial\Omega$ one
arrives at the definition of the space  
$B_p^{l-1/p} (\partial\Omega)$.

\smallskip

By the space of Sobolev multipliers $M(W^h_p (\Omega) \to W^l_p
(\Omega))$ we  mean the class of functions $\gamma$ such that $\gamma
u \in W_p^l (\Omega)$ for all $u \in W_p^h (\Omega)$. The space
$M(W^h_p (\Omega) \to W_p^l (\Omega)$) is endowed with the norm
\begin{equation}\label{v1}
\Vert \gamma  \Vert_{M(W_p^h(\Omega) \to W_p^l(\Omega))}= \sup \{\Vert\gamma
u  \Vert_{W_p^l(\Omega)}: \Vert u  \Vert_{W_p^h( \Omega)} \le 1\}.
\end{equation}
The notation $MW_p^l (\Omega)$ is used for $h =l$. Properties and applications of Sobolev multipliers are studied in \cite{[MS2]}.

\smallskip

We introduce the essential norm of a function $\gamma \in M (W_p^h
(\Omega) \to W_p^l (\Omega)$):
\begin{equation}\label{v2}
ess \Vert\gamma\Vert_{M(W^h_p( \Omega)\to W^l_p( \Omega))} =\inf_{\{T\}} \Vert\gamma -
T  \Vert_{W^h_p( \Omega) \to W^l_p( \Omega)},
\end{equation}
where $\{T\}$ is the set of all compact operators: $W^h_p (\Omega) \to
W^l_p (\Omega)$.

\smallskip

Analytic two-sided and one-sided estimates for the norm (\ref{v1}) and the essential norm (\ref{v2}) can be found in \cite{[MS2]}.

\subsection{Subclasses of Lipschitz graph domains}

Let $\Omega$ be a bounded Lipschitz graph domain. We introduce
the class $M^{l-1/p}_p$ $(l=2, 3, \dots)$ of  boundaries $\partial\Omega$,
satisfying the following condition.
For every point of  $\partial \Omega$ there exists an $n$-dimensional
neighborhood in which $\partial \Omega$ is specified (in a certain
Cartesian coordinate system) by a function $\varphi$ such that 
$$\nabla
\varphi \in MB^{l-1-1/p}_p(\mathbb{R}^{n-1}).$$
Furthermore, by
definition, $M^{1-1/p}_p$ is the class of bounded Lipschitz graph domains.

\smallskip

We say that $\partial\Omega$ belongs to the class $B^{l-1/p}_p$ if $\partial
\Omega$ can be locally specified by a function $\varphi\in
B^{l-1/p}_p(\mathbb{R}^{n-1})$. Since 
$$MB^{l-1-1/p}_p(\mathbb{R}^{n-1}) \subset B^{l-1-1/p}_{p,
\text{loc}}(\mathbb{R}^{n-1}),  \quad l \ge 2, $$
 and $C^{0,
1}(\mathbb{R}^{n-1}) \subset B^{l-1/p}_{p,
\text{loc}}(\mathbb{R}^{n-1})$,
it follows that  any bounded domain $\Omega$ with  $\partial\Omega\in M^{l-1/p}_p$ satisfies  $\partial\Omega\in B^{l-1/p}_p$.

\smallskip

According to Corollary  4.3.8 in \cite{[MS2]},  for $p(l-1)>n$ we have
\[ \| \nabla \varphi \|_{MB^{l-1-1/p}_p(\mathbb{R}^{n-1})}\sim \sup_{x
\in \mathbb{R}^{n-1}} \| \nabla \varphi  \|_{B^{l-1-1/p}_p({\cal B}_1(x))}\;,\]
where ${\cal B}_r(x)$ is a ball in $\Bbb{R}^{n-1}$ with radius $r$ and center $x$. 
Therefore,  the classes $M^{l-1/p}_p$ and $B^{l-1/p}_p$ coincide for
$p(l-1)>n$.

\smallskip

For a bounded Lipschitz graph domain $\Omega$,   by
$B^{l-l/p}_p(\partial \Omega)$ we denote the space of traces
 on $\partial \Omega$ of functions in $W^l_p(\Omega)$. Taking into
account an analogous fact for special Lipschitz domains of the class
$M^{l-1/p}_p$ (see Subsection 9.4.3 in \cite{[MS2]}), we obtain that 
$MB^{l-1/p}_p(\partial \Omega)$ is the space of traces of functions in
$MW^l_p(\Omega)$.

\smallskip

  In our subsequent  exposition the following additional condition on
$\Omega$  plays an important role.

 We say that $\partial\Omega$ belongs to the class $M^{l-1/p}_p(\delta)$ if  for every  point $O \in \partial
\Omega$ there exists a neighborhood $U$ and a special Lipschitz domain $G=\{z=(x, y):x \in \mathbb{R}^{n-1}, y > \varphi (x)\}$ such that
$U\cap \Omega=U \cap G$ and
\begin{equation}\label{def}
 \| \nabla \varphi \|_{MB^{l-1-1/p}_p( \mathbb{R}^{n-1})} \le
\delta\;.
\end{equation}

Obviously, the  boundaries in  $M^{l-1/p}_p(\delta)$ belong
to the class $M^{l-1/p}_p$ and, therefore, to the class $B^{l-1/p}_p$.

\smallskip

The following assertion gives a local characterization  of the class $M^{l-1/p}_p(\delta)$.
In its statement  we use the notion of the $(p,j)$-capacity of a
compact set $e$ in $\Bbb{R}^{n-1}$:
\begin{equation}\label{cap}
C_{p,j}(e) = \inf \{\Vert u \Vert^p_{W_p^j (\Bbb{R}^{n-1})} : u \in
C_0^\infty (\Bbb{R}^{n-1}), \,\,\, u \ge 1 \hbox{ on } e\}.
\end{equation} 

\noindent
For various properties of this capacity see the books by V. Maz'ya \cite{[Maz2]} and D.R. Adams and L.-I. Hedberg \cite{[AH]}.

\begin{proposition}\label{th13.18/1}
Let $p(l-1) \le n$. The class $M^{l-1/p}_p(\delta)$ admits the following equivalent 
description. For any point $O\in \partial \Omega$ 
there exists a neighborhood $U$ and a special Lipschitz domain 
$G=\{z=(x, y): x \in \mathbb{R}^{n-1}, y > \varphi(x)\}$ such that
$U \cap \Omega=U \cap G$ and 
\begin{equation}\label{eq13.18/2}
\lim_{\varepsilon \to 0} \Big( \sup_{e \subset {\cal B}_\varepsilon} \frac{\|
D_{l-1/p}(\varphi; {\cal B}_\varepsilon) \|_{L_p(e)}}{[C_{p, l-1-1/p}(e)]^{1/p}}
+ \| \nabla \varphi  \|_{L_\infty({\cal B}_\varepsilon)}\Big) \le c\,
\delta\;.
\end{equation}
Here ${\cal B}_\varepsilon$ is the ball with centre at $O$ and radius
$\varepsilon$,  $c$ is a constant which depends only on $l$, $p$, $n$, 
 and 
\[ D_{j-1/p} (\varphi; {\cal B}_\varepsilon )(x)= \Big( \int_{{\cal B}_\varepsilon} |
\nabla_{j-1} \varphi(x)-\nabla_{j-1} \varphi(y)|^p\frac{dy}{|x-y|^{n-2+p}}
\Big)^{1/p}\;.\]
\end{proposition}

Proposition \ref{th13.18/1} and properties  of the capacity
  lead to the following sufficient condition formulated in terms of the $(n-1)$-dimensional Lebesgue measure ${m}_{n-1}$. 

\begin{corollary}\label{coro13.18/1}
$(i)$ If $n > p(l-1)$ and 
\[ \lim_{\varepsilon \to 0} \Big( \sup_{e \subset {\cal B}_\varepsilon}
\frac{\| D_{l-1/p}(\varphi; {\cal B}_\varepsilon) 
\|_{L_p(e)}}{({m}_{n-1} e)^{[n-p(l-1)]/(n-1)p}} +\| \nabla
\varphi  \|_{L_\infty({\cal B}_\varepsilon)}\Big)< c \, \delta\;,\]
then $\partial\Omega\in M^{l-1/p}_p(\delta)$.

$(ii)$ If $n=p(l-1)$ and 
\[ \lim_{\varepsilon \to 0} \Big( \sup_{e \subset {\cal B}_\varepsilon} \|
D_{l-1/p} (\varphi; {\cal B}_\varepsilon)  \|_{L_p(e)} | \log
({m}_{n-1} e )|^{(p-1)/p}+\| \nabla \varphi 
\|_{L_\infty({\cal B}_\varepsilon)}\Big)< c\,\delta\;,\]
then $\partial\Omega\in M^{l-1/p}_p(\delta)$.
\end{corollary}

Now we present another test for the inclusion of a function into $M^{l-1/p}_p(\delta)$ involving the Besov space $B^m_{q, p}$.

We say that the boundary of a  Lipschitz graph domain  $\Omega$  belongs to
$B^{l-1/p}_{q, p}$ ($l=1, 2, \dots, )$ if, for any point of $\partial
\Omega$, there exists
a neighborhood in which $\partial \Omega$ is specified in Cartesian
coordinates by a function $\varphi$ satisfying
\[ \int_{\mathbb{R}^{n-1}} \Big( \int_{\mathbb{R}^{n-1}} |
\nabla_{l-1} \varphi(x+h)-\nabla_{l-1} \varphi(x)|^q\, dx \Big)^{p/q}\frac{dh}{
|h|^{2-n-p}}< \infty\;.\]

\begin{corollary}\label{coro13.18/2}
Let $p(l-1) \le n$ and let $\Omega$ be a bounded Lipschitz graph domain  with
 $\partial\Omega\in B^{l-1/p}_{q, p}$,  where 
$$q \in [p(n-1)/(p(l-1)-1), \infty] \quad 
{\rm if}  \,\,\, p(l-1) <n$$
 and 
 $$q \in (p, \infty] \qquad 
{\rm if} \,\,\, p(l-1)=n.$$
 Further, let $\partial \Omega$ be locally defined in
Cartesian coordinates by $y =\varphi(x)$, where $\varphi$ is a
function with a Lipschitz constant less than $c \, \delta$. Then 
 $\partial\Omega\in M^{l-1/p}_p(\delta)$.
\end{corollary}

Setting $q=\infty$ in Corollary \ref{coro13.18/2}, one obtains the  simple
sufficient condition for the inclusion into  $M^{l-1/p}_p(\delta)$ formulated in terms of the modulus 
of continuity $\omega_{l-1}(t)$ of $\nabla_{l-1} \varphi$:
\begin{equation}\label{eq13.18/3}
\int_0 \Bigl(\frac{\omega_{l-1}(t)}{t}\Bigr)^p\, dt< \infty\;.
\end{equation}
Since $B^{l-1/p}_{\infty, p} \subset B^{l-1/p}_p$, it follows that 
(\ref{eq13.18/3})  is sufficient for $\partial\Omega \in B^{l-1/p}_p$.

\subsection{General elliptic boundary value problem}

Consider either scalar or square matrix  differential operators

\begin{equation}\label{eq2t}
P(x,D_x) =  \sum_{\vert\alpha\vert\le 2 m} a_\alpha (x) D^\alpha_x,
\quad P_j (x, D_x) = \sum_{\vert\alpha\vert \le k_j } a_{\alpha j} (x)
D^\alpha_x,
\end{equation}
where $x \in \overline{\Omega}, \; 1 \le j \le m$, and $D_x = (i^{-1}\partial_1, \ldots , i^{-1}\partial_n)$.  We assume that the
operators $P$, $\text{tr }P_1, \dots, \text{tr }P_m$ form an elliptic
boundary value problem in every point
 $O \in \partial \Omega$ with respect to the hyperplane $y=0$ and that
$P$ is an elliptic operator in $\Omega$.  

\smallskip

The next  result is  proved essentially in the same manner as Theorem 14.3.1 in \cite{[MS2]}.

\begin{theorem}\label{TH1} 
$(i)$ Suppose that for any neighbourhood $U \subset \Bbb{R}^n$ there exist
operators
$$
P^U (D_z) = \sum_{\vert \alpha\vert = 2m} a^U_\alpha D_z, \quad P^U_j
(D_z) = \sum_{\vert\alpha\vert= k_j} a^U_{\alpha j} D_z
$$
with constant coefficients such that $\{P^U; {\rm Tr}\,  P_j^U\}$ is the
operator of an elliptic boundary value problem in the half-space $\{x
\in \Bbb{R}^{n-1}, y \ge 0\}$.

$(ii)$ Let 
$$
\sum_{\vert\alpha\vert = 2m} \Vert a_\alpha (z) - a^U_\alpha
\Vert_{L_\infty (U\cap \Omega)} \le \delta, 
$$
\begin{equation}\label{eq3t}
\sum_{\vert\alpha\vert \le 2m} ess \Vert
a_\alpha\Vert_{M(W_p^{l-\vert\alpha\vert}(\Omega) \to
W_p^{l-2m}(\Omega))} \le \delta,
\end{equation}
 
\noindent
where $l$ is integer, $l\ge 2m, 1<p<\infty$. The constant $\delta$ here and elsewhere is supposed to be small in
comparison with the norms of the inverse operators $\{P^U ; {\rm Tr}\, 
P_j^U\}^{-1}$ for all $U$. Further, let the
coefficients $a_{\alpha j}$ satisfy similar conditions with $2m$
replaced by $k_j$. 

$(iii)$ Let  the boundary $\partial\Omega$ belong to the class $M_p^{l-1/p}(\delta)$ if $p(l-1)\leq n$ or to the class $B_p^{l-1/p}$ if $p(l-1) >n$, in either case with $1<p<\infty$.

Then the operator 
\begin{equation}\label{eq5t}
\{P; {\rm Tr}\,  P_j\}: W^l_p (\Omega) \to W_p^{l-2m} (\Omega) \times
\prod^m_{j=1} B_p^{l-k_j-1/p} (\partial\Omega)
\end{equation}
is Fredholm. In particular, for all $u \in W^l_p(\Omega)$ the \emph{a priori} estimate
\begin{equation}\label{eq13.12/2}
\| u \|_{W^l_p(\Omega)} \le c \Bigl(\! \| Pu
\|_{W^{l-2m}_p(\Omega)}+\sum^m_{j=1} \|  {\rm Tr }\,  P_j u \|_{B^{l-k_j-1/p}_p(\partial 
\Omega)} + \| u\|_{L_1(\Omega)}\!\Bigr)
\end{equation}
holds; the last norm in the right-hand side  can be omitted in the case of a
unique solution.
\end{theorem}

\smallskip

Note that even the  roughest  sufficient condition  (\ref{eq13.18/3}) on the domain for $(iii)$  in Theorem \ref{TH1} to hold is sharp.
Let $\omega$ be  an increasing function
in $C[0, 1]$ such that $\omega(0) =0$,
$$\delta\int_\delta^1 \omega(t)\frac{dt}{t^2} + \int_0^\delta \omega(t)\frac{dt}{t} \leq c\, \omega(\delta),$$
and 
$$\int^1_0 \Bigl(\frac{\omega(t)}{t}\Bigr)^p\, dt =\infty.$$
It was
 shown in Section 4.4.3 of \cite{[MS2]}  that one can construct a function $\varphi$ on $\mathbb{R}^{n-1}$ such that

\smallskip

$(i)$ the continuity modulus   of $\nabla_{l-1} \varphi$ does not
exceed $c \,\omega$ with $c=\text{const}$;

\smallskip

$(ii)$ $\text{supp } \varphi \subset Q_{2 \pi}$, where $Q_d=\{ x \in
\mathbb{R}^{n-1}: |x_i | < d\}$;

\smallskip

$(iii)$ $\varphi \notin B^{l-1/p}_p(\mathbb{R}^{n-1})$.

\smallskip

Given  $\varphi$, one can construct  a bounded domain $\Omega$ in $\mathbb{R}^{n}$ such that the Neumann 
 problem
 \begin{equation}\label{Ne}
 \Delta v -v =g \quad {\rm in}\,\,  \Omega, \quad \partial v/\partial \nu =h \quad {\rm on}\,\, \partial\Omega
 \end{equation}
 with $g\in W_p^{l-2}(\Omega)$ and $h \in B_p^{l-1-1/p}(\partial\Omega)$
  may fail to be solvable  in $W^l_p(\Omega)$.

\smallskip

Next we describe conditions on the coefficients of (\ref{eq2t}) which are
equivalent to those formulated in Theorem \ref{TH1} and follow from the results
in \cite{[MS2]}, Ch. 7. In its formulation, we use the notion of capacity of a compact set in $\Bbb{R}^n$ defined similarly to (\ref{cap}). 

\begin{corollary}\label{COR4} Conditions $(\ref{eq3t})$ in Theorem $\ref{TH1}$ can be stated as follows:

$(i)$ The coefficients $a_\alpha$ with $\vert\alpha\vert = 2m$ are in the
class $W_p^{l -2m} (\Omega)$ if $p(l-2m) > n$ and satisfy the
inequality
$$
\Vert a_\alpha (z) - a_\alpha^U \Vert_{L_\infty(U \cap \Omega)} +
\lim_{\varepsilon \to 0} \sup_{\{e \subset \Omega : diam (e) \le
\varepsilon\}} {\Vert \nabla_{l-2m} a_\alpha  \Vert_{L_p(e)} \over 
[C_{p,l-2m} (e)]^{1/p}} \le \delta
$$
if $p (l - 2m) \le n$;

$(ii)$ The coefficients $a_\alpha$ with $\vert\alpha\vert < 2m$ are in the
class $W^{l-2m}_p (\Omega)$ if $p(l-\vert\alpha\vert) > n$ and
satisfy the inequality 
$$\lim_{\varepsilon \to 0}\Bigl( \sup_{\{e\subset \Omega : diam (e) \le
\varepsilon\}}
{\Vert \nabla_{l-2m} a_\alpha  \Vert_{L_p(e)} \over
[C_{p,l-\vert\alpha\vert} (e)]^{1/p}} +\sup_{x\in \Omega, \rho\leq \varepsilon} \rho^{2m-|\alpha| -n/p} \|a_\alpha\|_{L_p(B_\rho(x)\cap\Omega)}\Bigr)\leq \delta
$$
if $p(l - \vert\alpha\vert) \le n$. Here ${B}_\rho(x)$ is a ball in $\Bbb{R}^{n}$ with radius $\rho$ and center $x$. 
\end{corollary}

Various other sufficient conditions for (\ref{eq3t}) follow from the results in
\cite{[MS2]}.

\subsection{Dirichlet  problem in terms of traces}

Let us first consider the Dirichlet problem as a particular case of the general boundary value problem dealt with in Subsection 3.3. We write the scalar or square matrix elliptic  operator $P$ in the form
\begin{equation}\label{eq13.15.5/1}
Pu = \sum_{|\alpha|, |\beta| \le m} (-1)^{|\alpha|} D^\alpha
(A_{\alpha \beta}(z) D^\beta u)\;.
\end{equation}
Suppose that the  coefficients
$A_{\alpha \beta}$ are in   $C^{l-m}(\bar{\Omega})$, $l \ge m$,  and that 
 the G{\aa}rding
inequality
\begin{equation}\label{eq13.16.1/1}
\text{Re} \int_{\Omega} \sum_{|\alpha|=|\beta|=m} A_{\alpha \beta}
D^\alpha u \overline{D^\beta u}\, d z \ge c \, \| u 
\|^2_{W^m_2(\Omega)}
\end{equation}
holds for $u \in C^\infty_0(\Omega)$. 
 Assume that $\Omega$ is a Lipschitz graph domain.

\smallskip

We introduce a sufficiently small finite open covering $\{ U\}$ of
$\bar{\Omega}$ and a corresponding partition of unity $\{
\zeta_{U}\}$. Let 
$$P_{jU}=\partial^{j-1}/\partial y^{j-1}, \,\,\,  j=1,
\dots, m \quad {\rm  if}\,\,\, U \cap \partial \Omega \ne \varnothing$$
and 
$$P_{jU}=0\qquad {\rm  if} \,\,\, U \cap \partial \Omega =\varnothing.$$
 The Dirichlet
boundary conditions will be prescribed by the operators 
$$P_j=\sum_U
\zeta_U P_{jU}.$$

We give  a  formulation of the  Dirichlet problem. Let us 
  look for a function $u \in W^l_p(\Omega)$ such that
\begin{equation}\label{eq13.16.4/1}
Pu=f \quad \text{in} \quad \Omega, \qquad {\rm Tr }\,  P_j
u=f_j\quad\text{on}\quad \partial \Omega, \quad j=1, \dots, m\;,
\end{equation}
where $f$ and $f_j$ are  functions in  $W^{l-2m}_p(\Omega)$
and $B^{l+1-j-1/p}_p(\partial \Omega)$ respectively.

Here is the principal result relating the problem (\ref{eq13.16.4/1}) borrowed from Subsection 14.5.4 in \cite{[MS2]}.

\begin{theorem}\label{th13.16.4/1}
Let any of the following conditions hold:

$(\alpha)$ $m=1$, $p(l-1) \le n$; $\partial \Omega\in M^{l-1/p}_p(\delta)$;

$(\beta)$ $m=1$, $p(l-1)>n$; $\partial\Omega\in B^{l-1/p}_p$;

$(\gamma)$ $m>1$, $\partial\Omega \in M^{l-1/p}_p$ and $\partial \Omega$ is
locally defined by equations of the form $y=\varphi(x)$, where
$\varphi$ is a function with a small Lipschitz constant $($for 
$p(l-1)>n$, this is equivalent to $\partial\Omega \in B^{l-1/p}_p$$)$.

Then the operator
\[ \{P; {\rm Tr}\,  P_j\}: W^l_p(\Omega) \to W^{l-2m}_p(\Omega) \times
\prod^m_{j=1} B^{l+1-j-1/p}_p(\partial \Omega)\]
is an isomorphism.
\end{theorem}

 An example in  Subsection 14.6.1 of \cite{[MS2]}  shows that for $p(l-1) \le
n$ and for $m=1$ the condition  $\partial\Omega\in M^{l-1/p}_p(\delta)$ in part $(\alpha)$ of
Theorem \ref{th13.16.4/1} cannot be 
replaced by the assumption that $\partial\Omega$ belongs to the class
$M^{l-1/p}_p \cap C^{l-1}$. To be precise,   a 
domain $\Omega$ is constructed with $\partial\Omega \in M^{3/2}_2 \cap C^1$ for which the problem 
\begin{equation}\label{eq13.17.1/1}
-\Delta u=f\quad \text{in}\quad \Omega\;, \qquad {\rm Tr }\, u=0
\quad\text{on}\quad \partial \Omega
\end{equation}
is solvable in $W^2_2(\Omega)$  not for all $f \in L_2(\Omega)$. This
means that the smallness of  the seminorm $\| \nabla \varphi 
\|_{MB^{1/2}_2(\mathbb{R}^{n-1})}$ in the definition of  $M^{3/2}_2(\delta)$
is essential for the solvability of  problem (\ref{eq13.17.1/1}) in
$W^2_2(\Omega)$.

\smallskip

The next assertion, which  follows directly from the Implicit
Function Theorem 9.5.2 in \cite{[MS2]}, shows that the condition $\partial\Omega \in
B^{l-1/p}_p$ with $p(l-1)>n$ is necessary for the solvability of 
  problem (\ref{eq13.16.4/1}) in $W^l_p(\Omega)$ for the
operator
$P$ of higher than second order.

\begin{theorem}\label{th13.17.2/1}
Let $\Omega$ be a bounded Lipschitz domain  and let  $l$  be
integer, $l \ge 2m$, $p(l-1)>n$, $1<p< \infty$,  and $m>1$. If there
exists a solution $u \in W^l_p(\Omega)$  of the problem 
\begin{equation}\label{eq13.17.2/1}
Pu=0 \,\,\, {\rm in} \,\,\, \Omega, \,\,\, {\rm Tr} \, u=0, \quad
{\rm Tr} \, P_2 u=1, \,\,\, {\rm Tr} \, P_j u=0, \,\, j=3, \dots,
m\;,
\end{equation}
then $\partial\Omega \in B^{l-1/p}_p$.
\end{theorem}

Under the additional assumption $\partial\Omega \in C^{l-2, 1}$, 
the necessity of the inclusion $\partial\Omega \in B^{l-1/p}_p$ for $p(l-1) \le n$ is proved in Subsection 14.6.2 of \cite{[MS2]}.

\begin{theorem}\label{th13.17.2/2}
Let $\partial\Omega$ be in  the class $C^{l-2, 1}$ and let $l$  be
integer, $l \ge 2m$, $p(l-1) \leq n$, $1<p< \infty$,  and $m>1$. If 
there exists  a
solution  $u \in W^l_p(\Omega)$  of  problem
$(\ref{eq13.17.2/1})$, then $\partial\Omega \in B^{l-1/p}_p$.
\end{theorem}

A similar result for   the second order operator $P$ in Subsection 14.6.2 of \cite{[MS2]} runs as follows.

\begin{theorem}\label{th13.17.2/3}
Let $l$ be integer, $l \ge 2$, $1< p< \infty$, $m=1$, and $P1 \le
0$. Let $\Omega$ be a domain with $\partial\Omega\in C^1$ and let the normal to
$\partial \Omega$ satisfy the Dini condition. If, for a nonpositive
function $f \in C^\infty_0(\Omega)$, there exists a solution $u \in
W^l_p(\Omega)$ of the problem
\begin{equation}\label{eq13.17.2/2}
Pu=f \quad {\rm in} \quad \Omega,\quad {\rm Tr}\,  u=0\;,
\end{equation}
then  $\partial\Omega \in B^{l-1/p}_p$.
\end{theorem}

Note  that the convergence requirement (\ref{eq13.18/3}), the heaviest assumption on $\Omega$ made in Corollary \ref{coro13.18/2},  is in a sense a sharp  condition for
solvability of the Dirichlet problems (\ref{eq13.17.2/1}) 
and (\ref{eq13.17.2/2}) in $W^l_p(\Omega)$. The corresponding domain is constructed with the help of the same function $\omega$ as in the case of the Neumann problem (\ref{Ne}) (see Example 15.6.1 in \cite{[MS2]}).

\smallskip

Example 15.5.2 from \cite{[MS2]} shows that surfaces in the class $M^{l-1/p}_p(\delta)$ with $p(l-1)<n$ may have $${\rm conic} \,\,\, {\rm vertices}\,\,\, {\rm  if} \,\,\, n > p(l-1)$$
 and 
 $$s\mbox{-dimensional}\,\,\, {\rm  edges}\,\,\, {\rm if} \,\,\, s<n-p(l-1).$$

\smallskip

 Suppose that for any point $O\in \partial\Omega$ there exists a neighborhood $U$ such that $U\cap \Omega$ is $C^\infty$-diffeomorphic to the domain 
 $$\Bbb{R}^s\times \{(x,y): y>\varphi(x_{s+1}, \dots ,x_{n-1})\}, \qquad  0\leq s\leq n-2,$$
  i.e. the dimensions of boundary singularities are at most $n-1-s$. Then  (\ref{def}) is equivalent to 
$$\|\nabla \varphi  \|_{MB^{l-1-1/p}_p(\Bbb{R}^{n-1-s})}\leq c\, \delta$$

\noindent
and, in particular, it takes the form
$$\|\nabla \varphi  \|_{B^{l-1-1/p}_{p, \text{unif}}(\Bbb{R}^{n-1-s})}\leq c\, \delta, $$

\noindent
if $n-s<p(l-1)\leq n$. In other words, $\partial\Omega\in M^{l-1/p}_p(\delta)$ if and only if $(n-1-s)$-dimensional domain $\{(x,y): y>\varphi(x_{s+1}, \ldots ,x_{n-1})\}$ belongs to $M^{l-1/p}_p(c\,\delta)$.

\subsection{ Dirichlet
problem in a variational  formulation}

It turns out that for equations and systems of order higher than two, the formulation of the Dirichlet problem can be changed so that the solvability condition $\Omega\in M_p^{l-1/p}(\delta)$ is replaced by the better one $\Omega\in M_p^{l+1-m-1/p}(\delta)$. We comment on this referring to  Section 14.5 of \cite{[MS2]}.

\smallskip

Let $\Omega$ be open subset $\mathbb{R}^n$ and let $P$ be the operator
(\ref{eq13.15.5/1}), where $A_{\alpha \beta} \in
C^{l-m}(\bar{\Omega})$, $l \ge m$. Further,
let  the G{\aa}rding
inequality (\ref{eq13.16.1/1})
 hold for $u \in C^\infty_0(\Omega)$.

\smallskip

We say that $u \in W^l_p(\Omega)$ is a variational solution of the  Dirichlet
problem  if 
\begin{equation}\label{eq13.16.1/2}
Pu=f\;, \qquad u-g \in W^l_p(\Omega) \cap \ring{W}^m_p(\Omega)\;,
\end{equation}
where $f$ and $g$ are given functions in the spaces $W^{l-2m}_p(\Omega)$
and $W^l_p(\Omega)$ respectively.

\smallskip

By $W^{-k}_p(\Omega)$, $k=1, 2, \dots,$ we mean the space of linear
continuous functionals in $\ring{W}^k_{p'}(\Omega)$.

\smallskip

We present an a  priori estimate for  solutions of  
 problem (\ref{eq13.16.1/2}).

\begin{theorem}\label{th13.16.2/1}
If 

$(i)$ either $p(l-m) \le n$, $1 < p< \infty$ and  $\partial\Omega$ belongs to the class
 $M^{l+1-m-1/p}_p(\delta)$, or 
 
 $(ii)$ $p(l-m)>n$, $1 < p< \infty$, and $\partial\Omega \in B^{l+1-m-1/p}_p$, 
  
  \noindent
  then
\begin{equation}\label{eq13.16.2/1}
\| u  \|_{W^l_p( \Omega)} \le c\, (\|Pu\|_{W^{l-2m}_p( \Omega)} +\| u \|_{L_1( \Omega)})
\end{equation}
for all $ u \in (W^l_p \cap \ring{W}^m_p)(\Omega)$.
\end{theorem}

Next we state two corollaries of (\ref{eq13.16.2/1}).

\begin{proposition}\label{pro13.16.2/1}
Let $\Omega$ satisfy the conditions of  Theorem
$\ref{th13.16.2/1}$. 

$(i)$ If the kernel of the operator
\begin{equation}\label{eq13.16.2/2}
P:(W^l_p \cap \ring{W}^m_p)(\Omega) \to W^{l-2m}_p(\Omega)
\end{equation}
is trivial, then the norm $\| u \|_{L_1(\Omega)}$ in
$(\ref{eq13.16.2/1})$ can be omitted.

$(ii)$ The kernel of the operator $(\ref{eq13.16.2/2})$ is
finite-dimensional and the range of this operator is closed.
\end{proposition}

\begin{proposition}\label{pro13.16.2/2}
Let $\Omega$ satisfy the conditions of Theorem
$\ref{th13.16.2/1}$. Further, let $U$
and $V$ be open bounded subsets of $\mathbb{R}^n$, $\bar{U} \subset V$ and $u
\in (W^l_p \cap \ring{W}^m_p)(\Omega)$. Then
\[ \| u  \|_{W^l_p(U \cap \Omega)} \le c\, (\| Pu \|_{W^{l-2m}_p(V \cap
\Omega)} +\| u \|_{L_1(V \cap \Omega)} )\;.\]
\end{proposition}

Finally we give a theorem on the solvability of  (\ref{eq13.16.1/2}).

Let the G{\aa}rding inequality (\ref{eq13.16.1/1})
 hold for all $ u \in C^\infty_0(\Omega)$. Then,  as  is well known, 
the equation $P u=f$ with $ f \in W^{-m}_2(\Omega)$
is uniquely solvable in $\ring{W}^m_2(\Omega)$.  

\begin{theorem}\label{th13.16.3/1}
Let $\partial\Omega \in M^{l+1-m-1/p}_p$ for $p(l-m) \le n$
and let
$\partial\Omega$ belong to the class $B^{l+1-m-1/p}_p$ for $p(l-m) >n$.

$(i)$ If $f \in W^{l-2m}_p(\Omega) \cap W^{-m}_2(\Omega)$, $ g \in
W^l_p(\Omega) \cap W^m_2(\Omega)$, $1 < p< \infty$, and if $ u \in
W^m_2(\Omega)$
is such that $Pu=f$, $u-g \in \ring{W}^m_2(\Omega)$, then $u \in
W^l_p(\Omega)$ and $u-g \in \ring{W}^m_p(\Omega)$.

$(ii)$ The problem $(\ref{eq13.16.1/2})$ has one  and only one
solution $ u \in W^l_p(\Omega)$.
\end{theorem}

\subsection{Strong solvability of the Dirichlet problem for the  Stokes system}

Here we complement Section 2 by a solvability of the Dirichlet problem for the Stokes system in weighted Sobolev spaces of higher order.

\smallskip

 Let $l$ be noninteger, $l>1$. We use the notation $M_p^l(\delta)$ for the class of three-dimensional Lipschitz graph domains subject to
$$\|\nabla \varphi\|_{MB_{p}^{l-1}(\Bbb{R}^{2})} \leq \delta$$
for an arbitrary coordinate system on $\partial\Omega$, where $\delta$ is a positive number.

\smallskip

We conclude this section by stating strong solvability result for the Dirichlet problem for the Stokes system
\begin{eqnarray}\label{s1}
&&\Delta u -\nabla \pi = f, \quad {\rm div}\, u = g \quad {\rm in} \,\,\, \Omega\nonumber\\
\\
&& {\rm Tr}\, u =h \qquad\qquad {\rm on} \,\,\, \partial\Omega\nonumber.
\end{eqnarray}
We assume that $g$ and $h$ satisfy the compatibility condition (\ref{comp-CD}) and use the space $W_p^{m,a}$ introduced in Subsect. 1.3. The proof is essentially the same as that of Theorem 15.1.2 in \cite{[MS2]}.

\begin{theorem}\label{St}
Let $p\in (1, \infty)$, $a = 1-\{l\} -1/p$, where $l$ is noninteger, $l>1$. Suppose that $\partial\Omega\in B_{p}^l$ for $p(l-1)>2$ and $\partial\Omega\in M_p^l(\delta)$ with some $\delta = \delta (p,l)$ for $p(l-1)\leq 2$. Then, for every triple
$$(f,g,h) \in W_p^{[l]-1,a} (\Omega)\times W_p^{[l],a} (\Omega)\times B_{p}^l(\partial\Omega)$$
there exists a unique solution $(u,\pi)$ of the problem $(\ref{s1})$ in $W_p^{[l]+1,a} (\Omega)\times W_p^{[l],a} (\Omega)$. 
\end{theorem}

Note that under the conditions of the last theorem, the operator
\begin{eqnarray*}
&&W_p^{[l]+1,a} (\Omega)\times W_p^{[l],a} (\Omega)\ni  (u,\pi) \,\,\, \Longrightarrow\\
\\
&& (\Delta u -\nabla \pi, \, {\rm div}\, u, \, {\rm Tr}\, u) \in W_p^{[l]-1,a} (\Omega)\times W_p^{[l],a} (\Omega)\times B_{p}^l(\partial\Omega)
\end{eqnarray*}
is continuous.

\section{Asymptotic formulas for solutions to elliptic equations}

\subsection{Asymptotics of solutions near Lipschitz boundary}

Results of a new type were obtained by V. Kozlov and V. Maz'ya in \cite{[KM3]} for solutions of the Dirichlet problem for higher order elliptic equations in Lipschitz graph domains. We mean an explicit description of the asymptotic behaviour of solutions near a point ${\cal O}$ of the Lipschitz boundary. 

\smallskip

Consider the special Lipschitz graph domain
$$G= \{x=(x', x_n)\in {\mathbb{R}}^n: \, \, x_n >\varphi(x')\},$$
where $\varphi(0) =0$ and $\varphi$ has a small Lipschitz constant. The authors of \cite{[KM3]} study solutions of arbitrary strongly elliptic equation of order $2m$ with constant complex-valued coefficients
\begin{equation}\label{1}
L(\partial_x) u(x) = f(x) \qquad {\rm on} \,\, B_3\cap G
\end{equation}
with zero Dirichlet data on $(B_3\cap \partial G)\backslash {\cal O}$, where $B_\rho = \{x: \, |x|<\rho\}$ and $\partial _x = (\partial_{1}, \ldots , \partial_{n})$. It is assumed that the operator $L$ has no lower-order terms and the coefficient in front of $\partial_{n}^{2m}$ is equal to $(-1)^m$. 

\smallskip

One of the results in \cite{[KM3]} is the existence of a solution ${\cal U}\in W^{m}_2 (G)$ of the homogeneous equation (\ref{1})) which admits the asymptotic representation
\begin{eqnarray}\label{2}
&&{\cal U}(x) = \exp\Bigl( -\int_{|x|<|y'|<1} \varphi(y') \partial_{n}^m E(y',0) \, dy' + O\bigl(\int_{|x|}^1 \varkappa^2(\rho)\frac{d\rho}{\rho} \bigr)\Bigr)\nonumber\\
\\
&&\times \Bigl( (x_n -\varphi(x'))^m + O \bigl (|x|^{m+1-\varepsilon}\bigl(\int_{|x|}^1 \varkappa(\rho)\frac{d\rho}{\rho^{2-\varepsilon}} +1\bigr)\bigr)\Bigr).\nonumber
\end{eqnarray}
Here $\varepsilon$ is a poisitive constant,
$$\varkappa(\rho) = \sup _{|y'|<\rho} |\nabla \varphi (y')|,$$
and $E$ is the Poisson solution of the equation $L(\partial_x)E(x) =0$ in the upper half-space ${\mathbb{R}}^n_+$ which is positive homogeneous of degree $m-n$ and is subject to the Dirichlet conditions on the hyperplane $x_n =0$
$$\partial_{n}^j E =0, \quad {\rm} \,\, 0\leq j\leq m-2, \quad {\rm and}\,\,\,\, \partial_{n}^{m-1} E = \delta(x'),$$
where $\delta$ is the Dirac function.

\smallskip

It is also shown that a multiple of ${\cal U}$ is the main term in the asymptotic representation of an arbitrary solution $u$ if both $u$ and $f$ are subject to mild growth conditions near ${\cal O}$.  
\smallskip

Solutions with the weakest possible  singularity at ${\cal O}$ are studied in \cite{[KM3]} as well. The authors construct a solution ${\cal U}$ of the homogeneous equation (\ref{1}) which is subject to the asymptotic formula
\begin{eqnarray}\label{5}
&&{\cal U}(x) = \exp\Bigl( \int_{|x|<|y'|<1} \varphi(y') \partial_{n}^m E(y',0) \, dy' + O\bigl(\int_{|x|}^1 \varkappa^2(\rho)\frac{d\rho}{\rho} \bigr)\Bigr)\nonumber\\
\\
&&\times \Bigl( E(x_n -\varphi(x')) + O \bigl (|x|^{m-n+1-\varepsilon}\bigl(\int_{|x|}^1 \varkappa(\rho)\frac{d\rho}{\rho^{2-\varepsilon}} +1\bigr)\bigr)\Bigr).\nonumber
\end{eqnarray}

\smallskip

The asymptotic formulas (\ref{2}) and (\ref{5}) can be simplified under additional conditions on $\varkappa(\rho)$. Let, in particular, 
$$\int_{0}^1 \varkappa^2(\rho)\frac{d\rho}{\rho}  <\infty.$$
Then, in the special case of the polyharmonic equation 
$$(-\Delta )^m u =0 \qquad {\rm on} \,\,\, B_3\cap G,$$
any solution $u$ satisfying $|u(x)| =O (| x|^{m-n+1-\varepsilon})$ is subject to the following alternatives:
 either 
$$u(x) \sim C\, \frac{(x_n -\varphi(x'))^m}{|x|^n} \exp\Bigl( m\frac{\Gamma(n/2)}{\pi^{n/2}}\int_{|x|<|y'|<1} \varphi(y') \frac{dy'}{|y'|^n}\Bigr)$$
or 
$$u(x) \sim C\, {(x_n -\varphi(x'))^m} \exp\Bigl(- m\frac{\Gamma(n/2)}{\pi^{n/2}}\int_{|x|<|y'|<1} \varphi(y') \frac{dy'}{|y'|^n}\Bigr).$$

\subsection{Asymptotics of solutions to equations with discontinuous coefficients near a smooth boundary}

Proofs of the just mentioned results in \cite{[KM3]}  rely upon the papers \cite{[KM1]}, \cite{[KM2]} on the asymptotic formulas for solutions to the Dirichlet problem for arbitrary even order $2m$ strongly elliptic equations of divergence form near a point ${\cal O}$ at the smooth boundary. It is required only that the coefficients of the principal part of the operator have small oscillation near this point, and the coefficients in lower order terms are allowed to have singularities at the boundary.

\smallskip

We say a few words on the proof of asymptotic formulas in \cite{[KM1]}, \cite{[KM2]}. The elliptic equation  is transformed to a first-order evolution system with the matrix whose entries are partial differential operators on the hemisphere with time dependent coefficients. Thus, the question of asymptotics of solutions to the original Dirichlet problem is reduced to the study of the long-time behaviour of solutions of the evolution system just mentioned. The structure of the operator matrix in the system is rather complicated, because it has been obtained from a higher order partial differential equation in the variational form. Moreover, the study of this system  is aggreviated by the scantiness of information about the behaviour of the operator matrics at infinity. This difficulty is overcome by a right choice of function spaces, characterizing the solutions and the right-hand side of the evolution system by certain seminorms depending on time.  To obtain an asymptotic formula for the solution, the authors apply a particular spectral splitting of the system into one-dimensional and infinite-dimensional parts. A general asymptotic theory of differential equations with operator coefficients in Banach spaces which is the basis of \cite{[KM1]}, \cite{[KM2]}  is developed in \cite{[KM]}.

\smallskip

As an illustration, we describe  a corollary of the main result in \cite{[KM1]} concerning second order equations. Consider the uniformly elliptic equation
\begin{equation}\label{f1}
-{\rm div} \, (A(x)\, \nabla u(x)) = f(x) \qquad {\rm in} \,\, \Omega
\end{equation}
complemented by the Dirichlet condition
\begin{equation}\label{f2}
u=0\quad {\rm on} \,\, \partial \Omega,
\end{equation}
where $\Omega$ is a domain in $\Bbb{R}^n$ with smooth boundary. We assume that elements of the $n\times n$-matrix $A(x)$ are measurable and bounded complex-valued functions. One consideres a solution $u$ with a finite Dirichlet integral and suppose, for simplicity,  that $f=0$ in a certain $\delta$-neighborhood $\Omega_\delta = \{ x\in \Omega: \, |x|<\delta\}$ of the origin. Further, it is assumed that there exists a constant symmetric matrix $A$ with positive definite real part such that the function
$$\sigma (r) := \sup_{\Omega_r} \|A(x) -A\|$$
is sufficiently small for $r<\delta$. 

\smallskip

We introduce  the function
$${\cal R} (x) = \frac{\langle (A(x) -A)\nu, \nu\rangle -n \langle A^{-1}(A(x) -A)\nu, x\rangle \langle\nu, x\rangle \langle A^{-1}x, x\rangle^{-1}} {|S^{n-1}| ({\rm det} A)^{1/2} \langle A^{-1}x, x\rangle ^{n/2}},
$$
where  $|S^{n-1}|$ is the Lebesgue measure of the unit sphere in $\Bbb{R}^n$, $\langle z, \zeta\rangle = z_1 \zeta_1 + \ldots + z_n\zeta_n$ and $\nu$ is the interior unit normal at a point ${\cal O}$ on the boundary of $\Omega$. (For the notation $({\rm det} A)^{1/2}$ and $\langle A^{-1} x, x\rangle ^{n/2}$ see \cite{[H]}, Sect. 6.2).

\smallskip
The following asymptotic formula for an arbitrary solution of (\ref{f1}), (\ref{f2}) with finite energy integral is a special case of the main theorem in \cite{[KM1]}:
\begin{eqnarray}\label{f4}
&&u(x) = C \exp \Bigl(-  \int_{\Omega_\delta\backslash \Omega_{|x|}} {\cal R}(y) dy + O\bigl(\int_{|x|}^\delta \sigma(\rho)^2\frac{d\rho}{\rho}\bigr) \Bigr)\nonumber\\
&&\times \Bigl({\rm dist} (x, \partial \Omega) + O\bigl( |x|^{2-\varepsilon}\int_{|x|}^\delta \sigma(\rho)\frac{d\rho}{\rho^{2-\varepsilon}}\bigr)\Bigr) + O(|x|^{2-\varepsilon}),
\end{eqnarray}
where $C = const$ and $\varepsilon$ is a small positive number.

\smallskip

Using (\ref{f4}), it is an easy matter to derive sharp two-sided estimate for the H\"older exponent of $u$ at the origin. Another direct application of (\ref{f4}) is the following criterion. Under the condition
\begin{equation}\label{f5}
\int_0^\delta \sigma(\rho)^2\frac{d\rho}{\rho} <\infty
\end{equation}
all solutions $u$ are Lipschitz at the origin if and only if
\begin{equation}\label{f6}
\mathop{\hbox {{\rm lim} {\rm inf}}}_{r\to +0}\int_{\Omega_\delta\backslash \Omega_r} {\rm Re}\,  {\cal R} \, dx >-\infty.
\end{equation}
Needless to say, this new one-sided restriction (\ref{f6}) is weaker than the classical Dini condition at the origin. The complementary assumption (\ref{f5}) appeared previously in several papers dealing with other problems of the boundary behaviour of solutions to equation (\ref{f1}) (see the articles by E. Fabes, D.  Jerison,  and C. Kenig \cite{[FJK]}, by B.E. Dahlberg \cite{[Dah3]}, by C. Kenig \cite{[Ken]} {\it et al}).

\smallskip

Let $v$ be a solution of the equation (\ref{f1}) complemented by the Dirichlet condition
$$v=0 \quad {\rm on} \,\,\, \partial\Omega\backslash \{{\cal O}\}, \,\,\, {\cal O}\in \partial\Omega,$$
which has an infinite energy integral and the least possible singularity. We state a particular case of the main theorem in \cite{[KM1]} which is the following asymptotic representation for $v$:
\begin{eqnarray}\label{f44}
&&v(x) = C \exp \Bigl(  \int_{\Omega_\delta\backslash \Omega_{|x|}} {\cal R}(y) dy + O\bigl(\int_{|x|}^\delta \sigma(\rho)^2\frac{d\rho}{\rho}\bigr) \Bigr)\nonumber\\
&&\times \Bigl(\frac{{\rm dist} (x, \partial \Omega)}{\langle A^{-1}x, x\rangle^{n/2}} + O\bigl( |x|^{2-n-\varepsilon}\int_{|x|}^\delta \sigma(\rho)\frac{d\rho}{\rho^{2-\varepsilon}}\bigr)\Bigr) + O(|x|^{1-\varepsilon}),
\end{eqnarray}
where $C = const$ and $\varepsilon$ is a small positive number.

\smallskip

In general, theorems proved in \cite{[KM1]} and \cite{[KM2]} provide asymptotic formulas similar to (\ref{f4}) and (\ref{f44}) for solutions of the Dirichlet problem for the higher order uniformly elliptic equation with complex-valued coefficients
$$\sum_{0\leq |\alpha|, |\beta|\leq m} (-\partial_x)^\alpha ({\cal L}_{\alpha\beta} (x)\, \partial^\beta_x u(x) ) = f(x) \quad {\rm on} \,\,\, B^+_\delta,$$
where $B^+_\delta = \Bbb{R}^n_+ \cap B_\delta$, $\Bbb{R}^n_+  = \{x=(x', x_n)\in \Bbb{R}^n: x_n>0\}$ and $B_\delta=\{x\in \Bbb{R}^n: |x|<\delta\}$. The only a priori assumption on the coefficients ${\cal L}_{\alpha\beta}$ is smallness of the function
$$\sum_{|\alpha|=|\beta| =m} |{\cal L}_{\alpha\beta} (x) - L_{\alpha\beta}| + \sum_{|\alpha +\beta|<2m} x_n^{2m-|\alpha+\beta|} |{\cal L}_{\alpha\beta} (x)|,
$$
where $x\in B^+_\delta$ and $L_{\alpha\beta}$ are constants.
 
\subsection{Corollaries of the asymptotic formulas in Section 4.1}

The last section in \cite{[KM3]} concerns, in particular, solutions to the Dirichlet problem for elliptic equations of order $2m$ with constant coefficients in plane domains with a small Lipschitz constant of the boundary as well as arbitrary bounded plane convex domans $\Omega$. Let $\Omega$ be a bounded domain in ${\mathbb{R}}^2$.  
Consider a strongly elliptic operator with constant coefficients
$$L(\partial_x) = \sum_{0\leq k\leq 2m} a_k\, \partial_{1}^k\, \partial_{2}^{2m-k},$$
and denote by $w$ a weak solution to the Dirichlet problem
\begin{equation}\label{80}
L(\partial_x) w = f, \qquad w\in \ring W^{m}_2(\Omega).
\end{equation}
If $f\in W^{-m}_2(\Omega)$, this problem is uniquely solvable.  let us assume that 
$$\Omega\cap B_{2\delta_0} = \{(x_1, x_2): x_2 >\varphi(x_1), \,\,\, |x|<2\delta_0\},$$
where $\varphi$ is a Lipschitz function on $[-2\delta_0, \, -2\delta_0]$ and $\varphi (0) =0$. Note that one does not require the convexity of $\varphi$.

\smallskip

The next  result concerning solutions to  problem (\ref{80}), 
which stems  from (\ref{2}),  is as follows.

\begin{theorem}\label{Th3}
Suppose that the Lipschitz norm of $\varphi$ on $[-2\delta_0, \, -2\delta_0]$ does not exceed a certain constant depending on the coefficients of $L$. Let $f$ be equal to zero in $\Omega\cap B_{2\delta}$. Then, for all $\delta\in (0, \delta_0)$, $x\in \Omega\cap B_{\delta}$ and $k= 1, \ldots, m-1$, 
\begin{eqnarray}\label{81}
&&|\nabla_k w(x)| \leq c \, A(2\delta) |x|^{m-k}\\
&&\times \exp\Bigl( -a\int_{|x|}^\delta\frac{\varphi(\rho) - \varphi(-\rho)}{\rho^2} d\rho + b\int_{|x|}^\delta\max_{|t|<\rho} |\varphi'(t)|^2\frac{d\rho}{\rho}\Bigr).\nonumber
\end{eqnarray}
Here  
\begin{equation}\label{82}
A(\delta) = \delta^{-1-m} \|w\|_{L_2(\Omega\cap B_{\delta})}
\end{equation}
and the notation
$$a=\frac{1}{2\pi} \Im \sum_{1\leq k\leq m} (\zeta_k^+ - \zeta_k^-),$$
is used, where $\zeta_1^+, \ldots, \zeta_m^+$ and $\zeta_1^-, \ldots, \zeta_m^-$ are roots of the polinomial $L(1, \zeta)$ with positive and negative imaginary parts, respectively. This value of $a$ is best possible. By $b$ and $c$  positive constants depending only on $m$ and the coefficients of $L$ are denoted.
\end{theorem}

\smallskip

Note that for the operator $\Delta^m$ one has $\zeta_k^{\pm} = \pm i$, which implies $a=-m/\pi$. 

\medskip

The next assertion is a consequence of Theorem \ref{Th3} when the function $\varphi$ is convex.

\begin{theorem}\label{Th4}
Suppose that the function $\varphi$ describing the domain $\Omega$ near the point ${\cal O}$ is non-negative and convex, and $|\varphi'(\pm 2\delta)|$ does not exceed a sufficiently small constant $l_0$ depending on $m$ and  the coefficients of $L(\partial_x)$. Furthermore, let $f$ be zero in $\Omega\cap B_{2\delta}$ and let $w$ be a solution of $(\ref{80})$, which is extended by zero outside $\Omega$. Then
\begin{equation}\label{83}
\|\nabla_m w\|_{L_\infty(B_\delta)} \leq c\, A(8\delta),
\end{equation}
where $\delta <\delta_0/8$ and $A(\delta)$ is given by $(\ref{82})$.
\end{theorem}

One of the main results obtained in \cite{[KM3]} concerns the Green function $G_L$ of the Dirichlet problem for the operator $L$ with real coefficients. 

\begin{theorem}\label{Th5}
Let $\Omega$ be an arbitrary bounded convex domain in ${\mathbb{R}}^2$ and let the coefficients of $L$ be real. Then, for all $x$, $y$ in $\Omega$,
\begin{equation}\label{95}
\sum_{|\alpha| =|\beta|=m} |\partial_x^\alpha\, \partial_y^\beta G_L(x,y)|\leq C\, |x-y|^{-2},
\end{equation}
where $C$ is a positive constant depending on $\Omega$.
\end{theorem}

The case of complex coefficients is more complicated.

\begin{theorem}\label{th6}
Let $L$ be an arbitrary  strongly elliptic operator with complex coefficients. Suppose that $\Omega$ is a bounded convex domain such that the jumps of all angles between the exterior normal vector to $\partial\Omega$ and the $x$-axis be smaller than a constant depending on $m$ and the coefficients of $L(\partial_x)$. Then, for all $x$, $y$ in $\Omega$, estimate $(\ref{95})$ holds.
\end{theorem}

Theorem \ref{Th5} implies the following regularity result.

\begin{corollary}\label{cor6}
Let $\Omega$ be an arbitrary bounded convex domain in ${\mathbb{R}}^2$ and let the coefficients of $L$ be real. Then the solution $w$ of problem $(\ref{80})$ with $f\in W^{1-m}_q(\Omega)$, $q>2$, satisfies
\begin{equation}\label{100}
\sum_{|\alpha|\leq m} \| D^\alpha w\|_{L_\infty(\Omega)} \leq C\, \|f\|_{W^{1-m}_q(\Omega)},
\end{equation}
where the constant $C$ depends on $\Omega$, $m$, $q$, and the coefficients of  $L(\partial_x)$. 
\end{corollary}

Generally, this assertion does not hold for operators with complex coefficients. More precisely, if there exists an angle vertex on the boundary of a convex domain $\Omega$, one can construct a second order strongly elliptic operator $L(\partial_x)$ with complex coefficients such that the Dirichlet problem (\ref{80}) with $f\in C(\overline\Omega)$ has a solution with unbounded gradient (see \cite{[KMR]}, Sect. 8.4.3). By Theorem \ref{th6}, the statement of Corollary \ref{cor6} for $L$ with complex coefficients holds if the jumps of the normal vector are either absent or small.

\subsection{Classical asymptotics of solutions near a point of the domain}

Now, we present some results borrowed from Sect. 14.9 of the book by V. Kozlov and V. Maz'ya \cite{[KM]} 
 which are devoted to the asymptotic behaviour of solutions to elliptic equations near an interior point $O$ of the  domain. Here a modified Dini-type condition on the coefficients is introduced which guarantees the preservation of the asymptotics of solutions to the main part of the equation with coefficients frozen at $O$.

\smallskip

Let
$$P(D_x) = \sum_{|\alpha| =2m} p_\alpha\, D^\alpha_x$$
and let $G$ denote the Green matrix of this operator, i.e.  the solution of the system
\begin{equation}\label{89}
P(D_x) \, G(x) = I_l\, \delta(x) \qquad {\rm in} \,\, \Bbb{R}^n
\end{equation}
where $I_l$ is the $l\times l$ identity matrix and $\delta$ is the Dirac function. It is well-known (see F. John's book \cite{[Jo]}) 
 that $G$ admits the representation
$$G(x) = 
\begin{cases}
r^{2m-n}Q(\omega) \qquad\qquad\qquad\,\,\,  {\rm if} \,\,\, 2m\geq n, \,\, n\,\,\,{\rm odd}, \,\,\, {\rm or} \,\,\, 2m<n\\
R(x)\, \log r + r^{2m-n} S(\omega) \quad {\rm if} \,\,\, 2m\geq n, \,\,\, n\,\,\, {\rm even},
\end{cases}$$
where $Q$ and $S$ are smooth matrix-functions on the unit sphere in $\Bbb{R}^n$ and $R$ is a homogeneous polynomial matrix of degree $2m-n$.

\smallskip

Let us consider the elliptic operator
$$Q(x, D_x) = \sum_{|\alpha| \leq 2m} q_\alpha(x)\, D^\alpha_x$$
with measurable coefficients in the punctured ball $B_{r_0}\backslash \{0\}$. We introduce the function
$$S(r) = \sup_{K_r}\Bigl\{ \sum_{|\alpha| =2m} |q_\alpha(x) -p_\alpha | +  \sum_{|\alpha| <2m} |x|^{2m -|\alpha|} |q_\alpha(x)|\Bigr\},$$
where $K_r = \{x\in \Bbb{R}^n: e^{-1}r < |x| < r\}$, and assume that $S(r)$ does not exceed a small positive constant. We shall also use the notation
$$\|u\|_{{\cal W}_2^{2m}(K_r)} = \Bigl( \sum_{|\alpha| \leq 2m} r^{2|\alpha| -n} \, \|D^\alpha_x u\|^2_{L_2(K_r)}\Bigr)^{1/2}.$$ 

\smallskip

We formulate three theorems on the asymptotic representation  as $x\to 0$ for solutions $u\in W_{2, loc}^{2m}(B_{r_0}\backslash \{0\})$ of 
\begin{equation}\label{90}
Q(x, D_x) u =0 \quad {\rm on} \,\,\,  B_{r_0}\backslash \{0\}
\end{equation}
satisfying 
\begin{equation}\label{91}
\|u\|_{{\cal W}_2^{2m}(K_r)} = O (r^{k+\delta})
\end{equation}
with some $\delta >0$ and integer $k$. 

\begin{theorem}\label{9.1}
Let $2m <n$ and
\begin{equation}\label{92}
\int_0^{r_0} S(r)\, |\log r|^{\gamma -1} \frac{dr}{r} <\infty,
\end{equation}
where $\gamma$ is a positive integer.

$(i)$ If $k\geq 0$, then
\begin{equation}\label{93}
u(x) = \sum_{|\alpha|= k+1} c_\alpha\, x^\alpha + v(x),
\end{equation}
where $c_\alpha= const$ and
\begin{equation}\label{94}
\|v \|_{{\cal W}_2^{2m}(K_r)} = o (r^{k+1}  |\log r|^{1-\gamma}).
\end{equation}

$(ii)$ If $k\leq 2m-n-1$, then
\begin{equation}\label{95a}
u(x) = \sum_{|\alpha|=2m -n- k-1} C_\alpha\, D^\alpha_x G(x) + v(x),
\end{equation}
where $C_\alpha= const$ and  $G$ is the Green matrix introduced by $(\ref{89})$.

$(iii)$ If $k =  2m-n$, then
\begin{equation}\label{96}
u(x) = {\rm const} + v(x),
\end{equation}
where 
$$
\|v \|_{{\cal W}_2^{2m}(K_r)} = o( |\log r|^{1-\gamma}).
$$
\end{theorem}

\smallskip

 The asymptotics (\ref{96}) can be made more precise under the assumption that the operator $Q$ contains no derivatives of order $|\alpha| <s$, that is
$$Q(x, D_x) = \sum_{s\leq|\alpha| \leq 2m} q_\alpha(x)\, D^\alpha_x.$$
The formula (\ref{96}) can be replaced by
$$u(x) = \sum_{|\alpha|\leq s} c_\alpha\, x^\alpha + v(x),$$
where
$$
\|v \|_{{\cal W}_2^{2m}(K_r)} = o( r^s\, |\log r|^{1-\gamma}).
$$

\begin{theorem}\label{9.3}
Let $n$ be odd, $2m>n$ and let $S$ be subject to $(\ref{92})$. Then
\begin{equation}\label{97}
u(x) = \sum_{|\alpha|= k+1} c_\alpha\, x^\alpha + \sum_{|\beta |=2m -n- k-1} C_\beta\, D^\beta_x G(x) +
 v(x),
\end{equation}
where $c_\alpha$ and $C_\beta$ are constants and $v$ satisfies $(\ref{94})$. If either $k<-1$ or $k>2m-n$, then the first or the second sum in $(\ref{97})$ should be omitted.
\end{theorem}

\begin{theorem}\label{9.4}
Let $n$ be even, $2m>n$.

$(i)$ If $k\leq -2$ and $(\ref{92})$ holds, then $u$ satisfies $(\ref{95})$ with $v$ subject to $(\ref{94})$.

$(ii)$ If $k\geq 2m -n$ and $(\ref{92})$ holds, then $u$ satisfies $(\ref{93})$ with $v$ subject to $(\ref{94})$.

$(iii)$ Let $-1\leq k \leq 2m-n-1$ and let $(\ref{92})$ be valid with $\gamma\geq 2$. Then $u$ is represented by $(\ref{97})$ with
$$\|u\|_{{\cal W}_2^{2m}(K_r)} = o( r^{k+1}\, |\log r|^{3-\gamma}).
$$
\end{theorem}

\subsection{Asymptotics of solutions of the second order equation with square-Dini coefficients}

The asymptotic behaviour  of solutions near the  isolated point $O$ of the domain  was recently considered by V. Maz'ya and R. McOwen \cite{[MM1]}, \cite{[MM2]} for the case of the second order elliptic operator in nondevergence form
\begin{equation}\label{m1}
{\cal L}(x, D_x) u(x) = \sum_{1\leq i,j\leq n}a_{i j}(x) \, {\partial _i \partial _j\, u}.
\end{equation}
It is assumed that the coefficients have modulus of continuity $\omega$ satisfying the square-Dini condition
\begin{equation}\label{m2}
\int_0^1\omega^2(t)\frac{dt}{t} <\infty.
\end{equation}
If the coefficients are real, then, without loss of generality one can put $a_{ij}(0) =\delta_{i j}.$

\smallskip

An important role is played by the function
\begin{equation}\label{m3}
I(r) = \frac{1}{|S^{n-1}|} \int_{r<|z|<\varepsilon} \Bigl( {\rm trace}\,  ({\bf A}_z) - n\frac{\langle {\bf A}_z z, z\rangle}{|z|^2}\Bigr) \frac{dz}{|z|^n},
\end{equation}
where ${\bf A}_z$ stands for  the matrix $a_{ij}(z)$,  $\langle, \rangle$ is the inner product in $\Bbb{R}^n$,  and $\varepsilon$ is a sufficiently small positive number.

\smallskip

If the coefficients are subject to the usual Dini condition
\begin{equation}\label{m4}
\int_0^1\omega (t)\frac{dt}{t} <\infty,
\end{equation}
then, obviously, there exists a finite limit of $ I(r)$ as $r\to 0$, 
but (\ref{m4}) is not necessary for the existence of this limit.  In general, under (\ref{m2}), $I(r)$ may be unbounded as $r\to 0$, but, clearly, for every $\lambda>0$ there exists $C_\lambda$ such that 
\begin{equation}\label{m5}
|I(r)| \leq \lambda\, |\log r| + C_\lambda \quad {\rm for} \,\,\, 0<r<\varepsilon.
\end{equation}

\smallskip

The results in \cite{[MM1]} and \cite{[MM2]} are formulated in terms of the $L_p$-means
\begin{equation}\label{m6}
M_p(w,r) := \Bigl( \meanint_{r<|x|<2r} |w(x)|^p \, dx\Bigr)^{1/p}.
\end{equation}

\smallskip

The main theorem in \cite{[MM1]} can be stated  as follows.

\begin{theorem}\label{theo1.1}
Suppose that
$$|a_{ij}(x) -\delta_{i j}| \leq \omega(|x|) \quad {\rm as} \,\,\, x\to O,$$
where $\omega$ satisfies $(\ref{m2})$. For $p\in (1, \infty)$ and $\varepsilon >0$ sufficiently small, there exists a weak solution $Z\in L_{p,loc}(B_\varepsilon)$ of 
\begin{equation}\label{m7}
\sum_{1\leq i,j\leq n}{\partial _i \partial _j}\bigl( a_{i j} (x) \, Z \bigr) =0
\end{equation}
satisfying
\begin{equation}\label{m8}
Z(x) = e^{-I(|x|)} (1 + \zeta (x) ),
\end{equation}
where the remainder term $\zeta$ satisfies
\begin{equation}\label{m9}
M_p(\zeta, r) \leq c\, \max (\omega(r), \sigma(r)) \quad {\rm with} \,\,\, \sigma(r) := \int_0^r\frac{\omega^2(t)}{t} dt.
\end{equation}
Moreover, if $u\in  L_{p,loc}(\overline{B_\varepsilon}\backslash \{O\})$ is a weak solution of 
$$
 \sum_{1\leq i,j\leq n}{\partial _i \partial _j}\bigl( a_{i j} (x) \, u \bigr) =0
$$
in $B_\varepsilon$ subject to the growth condition 
$$M_p(u,r) \leq c\, r^{2-n+\varepsilon_0},$$
where $\varepsilon_0 >0$, then there exists a constant $C$, depending on $u$, such that
\begin{equation}\label{m10}
u(x) = C\, Z(x) +w(x),
\end{equation}
where the remainder term $w$ satisfies 
$$M_p(w,r) \leq c\, r^{1-\varepsilon_1}$$ for $0<r<\varepsilon$ and any $\varepsilon_1>0$.
\end{theorem}

Combining (\ref{m10}) and (\ref{m8}), one obtains the asymptotic representation
\begin{equation}\label{m11}
u(x) = c\, e^{-I(|x|)} |x|^n (1+o(1)) \quad {\rm as}\,\,\, x\to O.
 \end{equation}
 
 \smallskip
 
 Now let us turn to the behaviour of a solution to the homogeneous nondivergence equation with an isolated singularity at $x=O$. We assume that $n>2$. The results obtained in \cite{[MM2]} imply the existence of the solution 
 \begin{equation}\label{m12}
 Z(x) \sim \frac{|x|^{2-n}}{n-2}\, e^{I(|x|)} \quad {\rm as}\,\,\, x \to O.
 \end{equation} 
The behaviour of $I(|x|)$ as $x \to O$ not only controls the leading asymptotics of $Z(x)$ but whether we can solve the equation
$${\cal L}(x, D_x) F(x) =\delta(x).$$
There are three important cases to consider.

\smallskip

$$1. \quad  I(0) = \lim_{x\to O} I(|x|) \quad {\rm exists} \,\, {\rm and} \,\, {\rm is}\,\, {\rm finite}.$$
In this case $Z(x)$ may be scaled by a constant factor  to make it asymptotic to the fundamental solution for the Laplacian. In fact, the distributional equation
 \begin{equation}\label{m13}
 -{\cal L}(x, D_x) Z(x) = C\, \delta(x)
 \end{equation} 
 can be solved to find
 \begin{equation}\label{m14}
 C= |S^{n-1}|\, e^{I(0)}.
  \end{equation} 

\smallskip

$$2. \quad I(|x|) \to -\infty\quad {\rm as} \,\,\, x\to O.$$
We see that 
$$Z(x) = o(|x|^{2-n}) \quad {\rm as} \,\,\, x\to O,$$
and we can solve (\ref{m13}) to find $C=0$. Thus, in this case we obtain the interesting corollary 
that
 \begin{equation}\label{m15}
 {\cal L} u =0 \quad {\rm in} \,\,\, B_\varepsilon
 \end{equation}admits a solution $u=Z$ that is quite singular at $x=O$:
 $$|Z(x)| \geq C_\lambda |x|^{2-n+\lambda}$$
 for every $\lambda>0$. In particular, local regularity of solutions of the homogeneous equation (\ref{m15}) does not hold.
 
 \smallskip
 
 $$3.\quad I(|x|) \to \infty \quad {\rm as}\,\,\, x\to O.$$
 Now we find that
 $$Z(x) |x|^{n-2} \to \infty \quad {\rm as} \,\,\, x\to O,
 $$
 so this solution grows more rapidly than the fundamental solution for the Laplacian. Although $Z$ still satisfies (\ref{m13}), we can no longer find $C$.

\bibliographystyle{amsalpha}

\end{document}